\documentclass[11pt]{amsart}

\usepackage{comment}

\newtheorem{theorem}{Theorem}[section]

\newtheorem{corollary}[theorem]{Corollary}

\newtheorem{lemma}[theorem]{Lemma}

\newtheorem{proposition}[theorem]{Proposition}

\theoremstyle{definition}
\newtheorem{remark}[theorem]{Remark}
\newtheorem{question}[theorem]{Question}

\usepackage[colorlinks, bookmarks=true]{hyperref}
\usepackage{color,graphicx,shortvrb}
\usepackage{multicol}
\usepackage{enumerate}
\usepackage{amsmath}
\usepackage{amssymb}
\usepackage{tikz-cd}

\newcommand{\supp}{\mathrm{supp}}

\newcommand{\bob}{\boldsymbol{\beta}}

\newcommand{\K}[2]{\mathcal K_{{#1},{#2} } }

\newcommand{\N}{\mathbb{N}}

\newcommand{\R}{\mathbb{R}}

\tikzcdset{scale cd/.style={every label/.append style={scale=#1},
		cells={nodes={scale=#1}}}}

\newcommand{\triple}[1]{{\left\vert\kern-0.25ex\left\vert\kern-0.25ex\left\vert #1 
    \right\vert\kern-0.25ex\right\vert\kern-0.25ex\right\vert}}

\begin{document}
	
\title{Approximate ultrahomogeneity in $L_pL_q$ lattices}

\author{Mary Angelica Tursi}

\maketitle

\begin{abstract}
 	We show that for $1\leq p, q<\infty$ with $p/q\notin \N$, the doubly atomless separable $L_pL_q$ Banach lattice $L_p(L_q)$ is approximately ultrahomogeneous (AUH) over the class of its finitely generated sublattices.  The above is not true when $p/q \in \N$. However, for any $p\neq q$, $L_p(L_q)$ is AUH over the finitely generated lattices in the class $BL_pL_q$ of bands of $L_pL_q$ lattices.
\end{abstract}

\section{Introduction}\label{s:intro}

In this paper, we explore the homogeneity properties (or lack thereof) of the class of $L_pL_q$ lattices under various conditions. \\

The following is taken from \cite{henson07}: A Banach lattice $X$ is an \textbf{abstract $L_pL_q$ lattice}  if there is a measure space $(\Omega, \Sigma, \mu)$ such that $X$ can be equipped with an $L_\infty(\Omega)$-module and a map $N:X\rightarrow L_p(\Omega)_+$ such that
\begin{itemize}
	\item For all $\phi \in L_\infty(\Omega)_+$ and $x\in X_+$, $\phi \cdot x\geq 0$,
	\item  For all $\phi \in L_\infty(\Omega)$ and $x\in X$, $N[\phi\cdot x] = |\phi|N[x]$.
	\item  For all $x,y \in X$, $N[x+y] \leq N[x] +N[y] $
	\item If $x$ and $y$ are disjoint, $N[x+y]^q= N[x]^q +N[y]^q$, and if $|x| \leq |y|$, then $N[x]\leq N[y]$.
	\item  For all $x\in X$, $\|x\| = \|N[x]\|_{L_p}$.
\end{itemize}

When the abstract $L_pL_q$ space is separable, it has a concrete representation: Suppose $(\Omega, \Sigma, \mu)$ and $(\Omega', \Sigma', \mu')$ are measure spaces.  Denote by $L_p(\Omega; L_q(\Omega'))$ the space of Bochner-measurable functions $f:\Omega \rightarrow L_q(\Omega')$ such that the function $N[f]$, with $ N[f](\omega) = \| f(\omega) \|_q$ for $\omega \in \Omega$, is in $L_p(\Omega)$. The class of \textit{bands} in $L_pL_q$ lattices, which we denote by $BL_pL_q$, has certain analogous properties to those of $L_p$ spaces, particularly with respect to its isometric theory. \\

 $L_pL_q$ lattices (and their sublattices) have been extensively studied for their model theoretic properties in \cite{henson07} and \cite{henson11}.  It turns out that while abstract $L_pL_q$ lattices themselves are not axiomatizable, the larger class $BL_pL_q$ is axiomatizable with certain properties corresponding to those of $L_p$ spaces. For instance, it is known that the class of atomless $L_p$ lattices is separably categorical, meaning that there exists one unique atomless separable $L_p$ lattice up to lattice isometry. Correspondingly, the class of \textit{doubly atomless} $BL_pL_q$ lattices is also separably categorical; in particular, up to lattice isometry, $L_p([0,1]; L_q[0,1])$, which throughout will just be referred to as $L_p(L_q)$, is the unique separable doubly atomless $BL_pL_q$ lattice (see \cite[Proposition 2.6]{henson11}).   \\
 
 Additionally, when $p\neq q$, the lattice isometries of $L_pL_q$ lattices can be characterized in a manner echoing those of linear isometries over $L_p$ spaces (with $p \neq 2$).  Recall from \cite[Ch. 11 Theorem 5.1]{banach1932theorie} that a map $T: L_p(0,1) \rightarrow L_p(0,1)$ is a surjective linear isometry iff $T f (t) = h(t) f(\phi(t))$, where $\phi$ is a measure-preserving transformation and $h$ is related to $\phi$ through Radon-Nikodym derivatives.  If we want $T$ to be a \textit{lattice} isometry as well, then we also have $h$ positive (and the above characterization will also work for $p = 2$). In \cite{cambern1974isometries} (for the case of $q = 2$) and \cite{Sourour1977OnTI}, a corresponding characterization of linear isometries is found for spaces of the form $L_p(X; Y)$, for certain $p$ and Banach spaces $Y$.  In particular, for $L_pL_q$ lattices with $p\neq q$:  given $f\in L_p(\Omega; L_q(\Omega'))$, where $f$ is understood as a map from $\Omega$ to $L_q$, any surjective linear isometry $T$ is of the form
 
 $$ Tf(x) = S(x)\big(e(x) \phi f (x) \big),$$
 
 where $\phi$ is a set isomorphism (see \cite{cambern1974isometries} and \cite{Sourour1977OnTI} for definitions) $e$ is a measurable function related to $\phi$ via Radon-Nikodym derivatives, and $S$ is a Bochner-measurable function from $\Omega$ to the space of linear maps from $L_q$ to itself such that for each $x$, $S(x)$ is a linear isometry over $L_q$.\\
 
 In \cite{raynaud1986sous}, Raynaud obtained results on linear subspaces of $L_pL_q$ spaces, showing that for $1\leq q \leq p < \infty$, some $\ell_r$ linearly isomorphically embeds into $L_p(L_q)$ iff it embeds either to $L_p$ or to $L_q$.  However, when $1\leq p \leq q < \infty$, for $p \leq r \leq q$, the space $\ell_r$ isometrically embeds as a lattice in $L_p(L_q)$, and for any $p$-convex and $q$-concave Orlicz function $\phi$, the lattice $L_\phi$ embeds lattice isomorphically into $L_p(L_q)$.  Thus, unlike with $L_p$ lattices whose infinite dimensional sublattices are determined up to lattice isometry by the number of atoms, the sublattices of $L_pL_q$  are not so simply classifiable. \\
 
  In fact, the lattice isometry classes behave more like the $L_p$ linear isometries, at least along the positive cone, as is evident in certain equimeasurability results for $L_pL_q$ lattices.  In \cite{raynaud1986sous}, Raynaud also obtained the following on uniqueness of measures, a variation of a result which will be relevant in this paper: let $\alpha > 0, \alpha \notin \N$, and suppose two probability measures $\nu_1$ and $\nu_2$ on $\R_+$ are such that for all $s > 0$,  

\[ \int_0^\infty (t + s)^\alpha \ d\nu_1(t) = \int_0^\infty (t + s)^\alpha \ d\nu_1(t).  \]

Then $\nu_1 = \nu_2$.  Linde gives an alternate proof of this result in \cite{Linde1986UniquenessTF}. \\

Various versions and expansions of the above result appear in reference to $L_p$ spaces: for instance, an early result from Rudin generalizes the above to equality of integrals over $\R^n$: (\cite{rudin1976p}).  Assume that $\alpha > 0$ with $\alpha \notin 2\N$, and suppose that for all $\mathbf v\in R^n$, 

\[ \int_{\R^n} (1 + \mathbf v \cdot z)^\alpha \ d\nu_1(z) = \int_{\R^n} (1 + \mathbf v \cdot z)^\alpha \ d\nu_2(z) \]

Then $\nu_1 = \nu_2$.  An application of this result is a similar condition by which one can show that one collection of measurable functions $F:\R^n \rightarrow \R$, with $\mathbf f = (f_1,...,f_n)$ is equimeasurable with another collection $\mathbf g = (g_1,...,g_n)$ By defining $\nu_1$ and $\nu_2$ as pushforward measures of $F$ and $G$. In the case of $L_p$ spaces, if $f$ and $g$ are corresponding basic sequences whose pushforward measures satisfy the above for $\alpha = p$, then they generate isometric Banach spaces. Raynaud's result shows the converse is true for $\alpha \neq 4,6,8,...$.  A similar result in$L_p(L_q)$ from \cite{henson11} holds for $\alpha = p/q \notin \N$ under certain conditions, except instead of equimeasurable $\mathbf f$ and $\mathbf g$, when the $f_i$'s and $g_i's$ are mutually disjoint and positive and the map $ f_i \mapsto g_i$ generates a lattice isometry, $(N[f_1],...,N[f_n])$ and $(N[g_1],...,N[g_n])$ are equimeasurable.  \\

Recall that a space $X$ is \textit{approximately ultrahomogeneous }(AUH) over a class $\mathcal G$ of finitely generated spaces if for all appropriate embeddings $f_i;E\hookrightarrow X$ with $i = 1,2$, for all $E\in \mathcal G$ generated by $e_1,...,e_n \in E$, and for all $\varepsilon > 0$, there exists an automorphism $\phi:X\rightarrow X$ such that for each $1\leq j \leq n$,  $\| \phi\circ f_1(e_j) - f_2(e_j) \| < \varepsilon$. 
  
 \begin{center}
 	 
  \begin{tikzcd} 
  	X \arrow{rr}{\exists \phi} & \quad & X \\
  	& E \arrow{lu}{f_1} \arrow{ru}{f_2} & \quad
  \end{tikzcd}
  
 \end{center}

In the  Banach space setting, the embeddings are linear embeddings and the class of finitely generated spaces are finite dimensional spaces.   In the lattice setting, the appropriate maps are isometric lattice embeddings, and one can either choose finite dimensional or finitely generated lattices.  \\
   
The equimeasurability results described above can be used to show an approximate ultrahomogeneity of $L_p([0,1])$  over its finite dimensional linear subspaces only so long as $p\notin 2\N$ (see \cite{lusky1978some}).   Conversely, the cases where $p \in 2\N$ are not AUH over finite dimensional linear subspaces, with counterexamples showing linearly isometric spaces whose corresponding basis elements are not equimeasurability.   Alternate methods using continuous Fra\"iss\'e  Theory have since then been used to give alternate proofs of linear approximate ultrahomogeneity of $L_p$ for $p\notin 2\N$ (see \cite{fer20}) as well as lattice homogeneity of $L_p$ for all $1\leq p < \infty$ (see \cite{ben15}, \cite{fer20}).\\

This paper is structured as follows: in section \ref{s:prelims}, we first establish basic notation and give a characterization of finite dimensional $BL_pL_q$ lattices.  This characterization is used in subsequent sections for establishing both equimeasurability and ultrahomogeneity results later on. \\

In section \ref{s:blplq_homogeneity} we show that when $p\neq q$,  $L_p(L_q):= L_p(L_q)$  is AUH over the larger class of finite dimensional (and finitely generated) $BL_pL_q$ spaces.  This is done by characterizing representations of $BL_pL_q$ sublattices $L_p(L_q)$ in such a way that induces automorphisms over $L_p(L_q)$ making the homogeneity diagram commute.  The results here play a role in subsequent sections as well.\\

 In section \ref{s:equimeasurability}, we prove that if in addition $p/q\notin \N$, $L_p(L_q)$ is  also AUH over the class of its finitely generated sublattices.  First,  we determine the isometric structure of finite dimensional sublattices of $L_p(L_q)$ lattices by giving an alternate proof of \cite[Proposition 3.2]{henson11} showing that two sublattices $E$ and $F$ of $L_p(L_q)$, with the $e_i$'s and $f_i$'s each forming the basis of atoms, are lattice isometric iff $(N[e_1],...,N[e_n])$ and $(N[f_1),...,N[f_n])$ are equimeasurable.  The equimeasurability result allows us to reduce a homogeneity diagram involving a finite dimensional sublattice of $L_p(L_q)$ to one with a finite dimensional $BL_pL_q$ lattice, from which, in combination with the results in section \ref{s:blplq_homogeneity}, the main result follows. \\
 
  Section \ref{s:non_homgeneity} considers the case of $p/q \in \N$.  Here, we provide a counterexample to equimeasurability in the case that $p/q \in \N$ and use this counterexample to show that in such cases, $L_p(L_q)$ is not AUH over the class of its finite dimensional lattices.\\
  


\section{Preliminaries}\label{s:prelims}

We begin with some basic notation and definitions. Given a measurable set $A\subseteq \R ^n$, we let $\mathbf 1_A$ refer to the characteristic function over $A$.  For a lattice $X$, let $B(X)$ be the unit ball, and $S(X)$ be the unit sphere.\\

For elements $e_1,...,e_n$ in some lattice $X$, use bracket notation $<e_1,...,e_n>_L$ to refer to the Banach lattice generated by the elements $e_1,...,e_n$. In addition, we write $<e_1,...,e_n>$ without the $L$ subscript to denote that the generating elements $e_i$ are also mutually disjoint positive elements in the unit sphere. Throughout, we will also use boldface notation to designate a finite sequence of elements: for instance, for $x_1,...,x_n \in \R$ or $x_1,...,x_n \in X$ for some lattice $x$, let $\mathbf x = (x_1,...,x_n)$.  Use the same notation to denote a sequence of functions over corresponding elements: for example, let $(f_1,...,f_n) = \mathbf f$, or $(f_1(x_1),...f_n(x_n)) = \mathbf f(\mathbf x)$, or $(f(x_1),...,f(x_n)) = f(\mathbf x)$.  Finally, for any element $e$ or tuple $\mathbf e$ of elements in some lattice $X$, let $\bob(e)$ and $\bob(\mathbf e)$ be the band generated by $e$ and $\mathbf e$ in $X$, respectively. \\

Recall that Bochner integrable functions are the norm limits of simple functions $f:\Omega \rightarrow L_q (\Omega')$, with $f(\omega) = \sum_1^n  r_i \mathbf 1_{A_i} (\omega) \mathbf 1_{B_i}$, where $\mathbf 1_{A_i}$ and $\mathbf 1_{B_i}$ are the characteristic functions for $A_i \in \Sigma$ and $B_i \in \Sigma'$, respectively.  One can also consider $f\in L_p(\Omega; L_q(\Omega'))$ as a $\Sigma \otimes \Sigma'$-measurable function such that

\[ \|f \| = \bigg( \int_\Omega \|f(\omega) \|_q^p \ d \omega \bigg)^{1/p} = \bigg(\int_\Omega \bigg( \int_{\Omega'} | f(\omega, \omega') |^q \ d\omega' \bigg)^{p/q} \ d\omega \bigg)^{1/p} \]

Unlike the more familiar $L_p$ lattices, the class of abstract $L_pL_q$ lattices is not itself axiomatizable; however, the slightly more general class $BL_pL_q$ of bands in $L_p(L_q)$ lattices is axiomatizable. Additionally, if $X$ is a separable $BL_p L_q$ lattice, it is lattice isometric to a lattice of the form

\begin{align*} 
	& \bigg( \bigoplus_p L_p(\Omega_n; \ell_q^n)\bigg) \oplus_p L_p( \Omega_\infty; \ell_q) \\
	\oplus_p &  \bigg( \bigoplus_p L_p(\Omega'_n; L_q(0,1) \oplus_q \ell_q^n) \bigg) \\
	\oplus_p & \ L_p( \Omega'_\infty; L_q(0,1) \oplus_q\ell_q).
\end{align*}

$BL_pL_q$ lattices may also contain what are called \textit{base disjoint} elements.  $x$ and $y$ are base disjoint if $N[x] \perp N[y]$.  Based on this, we call $x$ a \textit{base atom} if whenever $0 \leq  y , z \leq x$ with $y$ and $z$ base disjoint, then either $N[y] = 0 $ or $N[z] = 0$.  Observe this implies that $N[x]$ is an atom in $L_p$. Alternatively, we call $x$ a \textit{fiber atom} if any disjoint $0 \leq y,z \leq x$ are also base disjoint. Finally, we say that $X$ is \textit{doubly atomless} if it contains neither  base atoms nor fiber atoms.\\ 
  
 Another representation of $BL_pL_q$ involves its finite dimensional subspaces. We say that $X$ is an $(\mathcal L_p\mathcal L_q )_\lambda$ lattice, with $\lambda \geq 1$ if for all disjoint $x_1,...,x_n \in X$ and $\varepsilon > 0$,  there is a finite dimensional $F$ of $X$ that is $(1+\varepsilon)$-isometric to a finite dimensional $BL_pL_q$ space containing $x'_1,...,x'_n$ such that for each $1\leq i\leq n$, $\|x_i - x'_i \| < \varepsilon$. Henson and Raynaud proved that in fact, any lattice $X$ is a $BL_pL_q$ space  iff $X$ is $(\mathcal L_p \mathcal L_q)_1$ (see \cite{henson07}). This equivalence can be used to show the following:

\begin{proposition}(Henson, Raynaud)\label{p:finite_representation2}
	If $X$ is a separable $B L_p L_q$ lattice, then it is the inductive limit of finite dimensional $BL_pL_q$ lattices.
\end{proposition}

The latter statement is not explicitly in the statement of Lemma $3.5$ in \cite{henson07}, but the proof showing that any $B L_p L_q$ lattice is $(\mathcal L_p \mathcal L_q)_1$ was demonstrated by proving the statement itself.  \\

Throughout this paper, we refer to this class of finite dimensional $BL_pL_q$ lattices as $B\K{p}{q}$.  Observe that if $E\in B\K{p}{q}$,  then it is of the form $\oplus_p (\ell_q^{m_i})_1^N$ where for $1\leq k \leq N$, the atoms $e(1,1),..., e(k,m_k)$ generate $\ell_q^{m_k}$. 

\begin{proposition}\label{p:BLpLq-representation}
	Let $E$ be a $B\K{p}{q}$ sublattice of $L_p(L_q)$ with atoms $e(k,j)$ as described above.  Then the following are true:
	
	\begin{enumerate}
		\item There exist disjoint measurable $A(k) \subseteq [0,1]$ such that for all $i$, $\supp(e(k,j)) \subseteq A(k) \times [0,1]$, 
		
		\item For all $k$ and for all $j,j'$, $N[e(k,j)]= N[e(k,j')]$.  
	\end{enumerate}
	Conversely, if $E$ is a finite dimensional sublattice of $L_p(L_q)$ satisfying properties (1) and (2), then $E$ is in $B\K{p}{q}$.
	
\end{proposition} 

In order to prove this theorem, we first need the following lemma:

\begin{lemma}\label{l:positive_p_disjoint}
	Let $0 < r < \infty$, with $r \neq 1$.  suppose $x_1,...,x_n \in L_r+$ are such that 
	
	\[ \| \sum_1^n  x_k \|_r^r = \sum \|x_k\|^r_r  \]
	
	Then the $x_i$'s are mutually disjoint.
\end{lemma}

\begin{proof}
	If $r < 1$, then	
	 \begin{align} 
	\int x_i(t)^r + x_j(t)^r \ dt =  \|x_i\|_r^r + \|x_j\|_r^r = 
	   \int (x_i(t) + x_j(t))^r \ dt 
	\end{align}
Now observe that for all $t$,  $(x_i(t) + x_j(t) )^r \leq x_i(t)^r + x_j(t)^r$, with equality iff either $x_i(t) = 0$ or $x_j(t) =  0$, so $(x_i + x_j)^r - x_i^r - x_j^r \in L_1+$.  Combined with the above equality in line (1), since $\|(x_i + x_j)^r - x_i^r - x_j^r \|_1 = 0, $ it follows that $x_i(t)^r + x_j(t)^r   = (x_i(t) + x_j(t))^r$ a.e., so $x_i$ must be disjoint from $x_j$ when $i \neq j$.  \\

If $r> 1$,  proceed as in the proof for $r < 1$, but with the inequalities reversed, given that in this instance $x_i(t)^r + x_j(t)^r \leq (x_i(t) +x_j(t))^r$ for all $t$.

\end{proof}

\begin{remark}
	The above implies that a $BL_pL_q$ lattice $X$ is base atomless if it contains no bands lattice isometric to some $L_p$  or $L_q$ space.  Indeed, if there were a base atom $e$, then any two $0 \leq x \perp y \leq e$ would have to have $N$-norms multiple to each other, so $<x,y>$ is lattice isometric to $\ell_q^2$.  Resultantly, the band generated by $e$ is an $L_q$ space. Similarly, if $e$ is a fiber atom, then any  $0 \leq x\perp y \leq e$ is also base disjoint, which implies that the band generated by $e$ is an $L_p$ space.\\
\end{remark}

We now conclude with the proof of Proposition \ref{p:BLpLq-representation}:
\begin{proof}[Proof of Proposition \ref{p:BLpLq-representation}]
	 
	Observe that for each appropriate pair $(k,j)$, 
	\[ \bigg( \int_0^1 N[e(k,j)]^{p}(s) \ ds \bigg)^{q/p} = \| N^q[e(k,j)] \|_{p/q} = 1   \]
	  For notational ease, let $E(k,j) = N^q[e(k,j)]$. Pick $j_1,...,j_n$ with each $j_k \leq m_k$. Then, by disjointness of the $e(k,j)$'s,  for all $(a_k)_k \geq 0$ and all $x = \sum_k a_k e(k,j_k)$,  
	  \begin{align*}
	  \| \sum a_k e(k,j_k) \|^q & =  \bigg( \int_0^1 \bigg(  \sum_k  a_k^q  E(k,j_k)(s) \bigg)^{p/q} \ ds  \bigg)^{q/p}\\
	   & = \bigg|\bigg| \sum a_k^q E(k,j_k)  \bigg|\bigg|_{p/q} .
	  \end{align*}
	  Now since the $e(k,j_k)$'s are isometric to $\ell_p$,
	  \begin{align*}
		\bigg|\bigg| \sum a_k^q E(k,j_k)  \bigg|\bigg|_{p/q}^{p/q}  = \sum_ i a_k^p = \sum_k (a_k^{q})^{p/q} = \sum_k \| a_k^q E(k,j_k) \|_{p/q}^{p/q}.
	  \end{align*}
	  Since the $E(k,j_k)$'s are all positive and $p\neq q$, by Lemma \ref{l:positive_p_disjoint}, the $E(k,j_k)$'s are disjoint, that is, the $e(k,j_k)'s$ are base disjoint. \\
	  
	  For $1\leq k\leq N$, let $A(1),...,A(n)$ be mutually disjoint measurable sets each supporting each $E(k,j)$ for $1 \leq j \leq n_k$.  Then each $e(k,j)$ is supported by $A(k)\times [0,1]$.  Now we prove (2).  Fix $k$,  Then using similar computations as above, and since the $e(k,j)$'s for fixed $k$ generate $\ell_q^{m_k}$: 
	  
	  \begin{align*}
	  \| \sum_j a_j e(k,j)  \|^q & = \bigg|\bigg| \sum_j a^q_j E(k,j) \bigg|\bigg|_{p/q} = \sum_j a_j^q = \sum_j a_j^q \|E(k,j)\|_{p/q}
	  \end{align*}
	  
	  By Minkowski's inequality, as $p\neq q$, equality occurs only when $E(k,j)(s)= E(k,j')(s)$ a.e. for all $1\leq j,j' \leq n_i$. \\

	  To show the converse, it is enough to give the computation: 
	  
	  \begin{align*}
	  \| \sum_{k,j} a(k,j) e(k,j) \| & = \bigg( \int_0^1 \bigg[ \int \bigg(\sum_{k,j} a(k,j)e(k,j)(s,t) \bigg)^q \ dt \bigg]^{p/q} \ ds \bigg)^{1/p} \\
	  & = \bigg( \sum_k \int_0^1 \bigg[  \sum_{j=1}^{n_i} |a(k,j)|^q E(k,j)(s) \bigg]^{p/q} \ ds \bigg)^{1/p} \\
	  & = \bigg( \sum_k \bigg[  \sum_{j=1}^{n_k} |a(k,j)|^q  \bigg]^{p/q} \int_0^1  E(k,1)^{p/q}(s) \ ds \bigg)^{1/p} \\
	  & = \bigg( \sum_k \bigg[  \sum_{j=1}^{n_k} |a(k,j)|^q  \bigg]^{p/q} \bigg)^{1/p} 
	  \end{align*}	  
\end{proof}

The following results will allow us to reduce homogeneity diagrams to those in which the atoms $e(k,j)$ of some $E\in B\K{p}{q}$ are mapped by both embeddings to characteristic functions of measurable $A(k,j) \subseteq [0,1]^2$. In fact, we can further simplify such diagrams to cases where $E$ is generated by such $e(k,j)$'s which additionally are \textit{base-simple}, i.e., $N[e(k,j)]$ is a simple function.

\begin{proposition}\label{p:homogeneity-one-element}
	Let $1 \leq p \neq q <\infty$ and let $e\in S(L_p(L_q))_+$ be an element with full support over $[0,1]^2$.  Then there exists a lattice automorphism $\phi$ from $L_p(L_q)$ to itself such that $\phi(\mathbf 1) = e$.  Furthermore, $\phi$ can be constructed to bijectively map both simple functions to simple functions and base-simple functions to base-simple functions.  
\end{proposition}

\begin{proof}
	
	The proof is an expansion of the technique used in Lemma 3.3 from \cite{fer20}.  Given a function $g(y) \in {L_q}_+$, define $\tilde g(y)_q$ by $\tilde g(y)_q = \int_0^y g(t)^q \ dt$, and for notation, use $e_x(y) = e(x,y)$. Since $e$ has full support, we may assume that for all $0 \leq x \leq 1$, $N[e](x) > 0$. From there, Define $\phi$ by
	
	\[ \phi(f)(x,y) = f \bigg(\widetilde{ N[e]}(x)_p, \frac{\tilde e_x(y)_q}{N^q[e](x)} \bigg) e(x,y) \]
	
	$e \geq 0$ and the rest of the function definition is a composition, so $\phi$ is a lattice homomorphism.  To show it is also an isometry, simply compute the norm, using substitution in the appropriate places:
	
	\begin{align*}
		\| \phi(f) \|^p = & \int_0^1 \bigg| \int_0^1 f \bigg(\widetilde{ N[e]}(x)_p, \frac{\tilde e_x(y)_q}{N^q[e](x)} \bigg)^q e(x,y)^q \ dy \bigg|^{p/q} \ dx \\
		= &  \int_0^1 \bigg| \int_0^1 f (\widetilde {N[e]}(x)_p, y)^q  \ dy \bigg|^{p/q} N^p[e](x) \ dx \\
		= & \int_0^1 N[f](\widetilde{ N[e]}(x)_p)^{p} N^p[e](x) \ dx \\
		= & \int_0^1 N^p[f](x) \ dx = \|f\|^p.
	\end{align*}
	
	To show surjectivity, let $B\subseteq [0,1]^2$ be a measurable set.   Note that any $(x',y')\in [0,1]^2$ can be expressed as $(\widetilde{N[e]}(x)_p, \frac{\tilde e_x(y)_q}{N^q[e](x)}) $ for some $x,y$, since $\widetilde{N[e]}(x)_p$ is an increasing continuous function from 0 to 1, while $\tilde e_x(y)_q$ is continuously increasing from 0 to $N^q[e](x)$.  Thus there exists $B'$ such that $\phi(\mathbf 1_{B'}) = \mathbf 1_B\cdot e$, implying that $\phi$'s image is dense in the band generated by $\bob(e) = L_p(L_q)$ since $e$ has full support. Therefore, $\phi$ is also surjective. \\
	
	Finally, $\phi$ consists of function composition into $f$ multiplied by $e$, so if $e$ and $f$ are simple, then it has a finite image, so if $f$ is simple, then the product is also simple,  $\phi$ maps simple functions to simple functions,  Conversely, if $\phi(f)$ is simple, then $\phi(f)/e$ is also simple.  Thus $f \bigg(\widetilde {N[e]}(x)_p, \frac{\tilde e_x(y)_q}{N[e](x)} \bigg)$ has a finite image. It follows that $f$ itself has a finite image.  \\
	
	Using similar reasoning, if $N[e]$ is simple, then whenever $N[f]$ is simple,  $N[\phi(f)]$ must also be simple, and likewise the converse is true, since by the computation above, $N[\phi(f)](x) = N[f](\widetilde{N[e]}(x)_p) \cdot N[e](x)$.
\end{proof}

\section{Approximate Ultrahomogeneity of $L_p(L_q)$ over $BL_pL_q$ spaces}\label{s:blplq_homogeneity}
In this section, we show that for any $1\leq p \neq q <\infty$, $L_p(L_q)$ is AUH over $B\K{p}{q}$.\\

Let $\mathbf f:= (f_1,...,f_n)$ and $\mathbf g:= (g_1,...,g_n)$ be sequences of measurable functions and let $\lambda$ be a measure in $\R$. Then we say that $\mathbf f$ and $\mathbf g$ are \textit{equimeasurable} if for all $\lambda$-measurable $B\subseteq \R^n$, 

\[ \lambda(t: \mathbf f(t) \in B) = \lambda(t: \mathbf g(t) \in B)  \]

We also say that functions $\mathbf{f}$ and $\mathbf g$ in $L_p(L_q)$ are \textit{base-equimeasurable} if $N(\mathbf f)$ and $N(\mathbf g)$ are equimeasurable. \\

Lusky's main proof in \cite{lusky1978some} of linear approximate ultrahomogeneity in $L_p(0,1)$ for $p \neq 4,6,8,...$ hinges on the equimeasurability of generating elements for two copies of some $E = <e_1,...,e_n>$ in $L_p$ containing $\mathbf 1$.  But when $p = 4,6,8,...$, there exist finite dimensional $E$ such that two linearly isometric copies of $E$ in $L_p$ do not have equimeasurable corresponding basis elements.  However, if homogeneity properties are limited to $E$ with mutually disjoint basis elements, then $E$ is linearly isometric to $\ell_p^n$, and for all $1\leq p <\infty$, $L_p$ is AUH over all $\ell_p^n$ spaces.  Note that here, an equimeasurability principle (albeit a trivial one) also applies:  Any two copies of $\ell_p^n = <e_1,..., e_n>$ into $L_p(0,1)$ with $\sum_k e_k = n^{1/p}\cdot \mathbf 1$  have (trivially) equimeasurable corresponding basis elements to each other as well. \\

In the $L_p(L_q)$ setting, similar results arise, except rather than comparing corresponding basis elements $f_i(e_1),...,f_i(e_n)$ of isometric copies $f_i(E)$ of $E$, equimeasurability results hold in the $L_q$-norms $N[f_i(e_j)]$ under similar conditions, with finite dimensional $BL_pL_q$ lattices taking on a role like $\ell_p^n$ does in $L_p$ spaces. \\

The following shows that equimeasurability plays a strong role in the approximate ultrahomogeneity of $L_p(L_q)$ by showing that any automorphism fixing $\mathbf 1$ preserves base-equimeasurability for characteristic functions:

\begin{proposition}\label{p:1-1-isometry_representation}
	Suppose $p \neq q$,  and let $T:L_p(L_q)$ be a lattice automorphism with $T(\mathbf 1) = \mathbf 1$.  Then there exists a function $\phi \in L_p(L_q)$ and a measure preserving transformation $\psi$ over $L_p$ such that for a.e. $x\in [0,1]$, $\phi(x,\cdot)$ is also a measure preserving transformation inducing an isometry over $L_q$, and for all $f$, 
	\[ Tf(x,y) = f(\psi(x), \phi(x,y)). \]
	Furthermore, for all measurable $B_1,...,B_n \subseteq [0,1]^2$ with $\mathbf 1_{B_i}$'s mutually disjoint, $(\mathbf 1_{B_1},..., \mathbf 1_{B_n})$ and $(T\mathbf 1_{B_1},..., T\mathbf 1_{B_n})$ are base-equimeasurable.
\end{proposition}

\begin{proof}
	By the main result in \cite{Sourour1977OnTI}, there exists a strongly measurable function $\Phi:[0,1]\rightarrow B(L_q)$, a set isomorphism $\Psi$ over $L_p$ (see \cite{Sourour1977OnTI} for a definition on set isomorphisms), and some $e(x) \in L_p$ related to the radon-Nikodym derivative of $\Psi$ such that
	
	\[ Tf(x) (y) = \Phi(x)(e(x)\Psi f(x)) (y), \] 
	
	and for a.e. $x$, $\Phi(x)$ is a linear isometry over $L_q$. Observe first that $T$ sends any characteristic function $1_{A\times [0,1]} \in L_p(L_q)$ constant over $y$ to characteristic function $\mathbf 1_{\psi(A)\times [0,1]}$ for some $\psi(A) \subseteq [0,1]$, so since $1_{A\times [0,1]} \in L_p(L_q)$ is constant over $y$, we can just refer to it as $\mathbf 1_A$. Also, since $T$ is a lattice isometry, $\mu(A) = \mu(\psi(A))$, so $\psi$ is measure preserving. Finally, observe that $N[\mathbf 1_{A}] = \mathbf 1_A$. Thus, for any simple function $g:= \sum c_i\mathbf 1_{A_i} \in L_p(L_q)_+$ constant over $y$ with the $A_i$'s mutually disjoint, we have $N[g] = g$, and $Tg = g'$. Then for all $x$, $$N[g'](x) = N[Tg](x) = N[\Phi(x)(e g')](x) = e(x) N[\Phi(x)(g')][x] = |e(x)|N[g'](x)$$ 
	 
	 It follows that $|e(x)| = 1$. We can thus adjust $\Phi$ by multiplying by $-1$ where $e(x) = -1$.  Note also that $\Phi$ acts as a lattice isometry over $L_p$ when restricted to elements constant over $y$, so by Banach's theorem in \cite{banach1932theorie}, the map  $\Phi f (x)$ can be interpreted as $\Phi(x)(\ f(\psi(x))\ )$, where $\psi$ is a measure preserving transformation over $[0,1]$ inducing $\Psi$.  By Banach's theorem again for $\Phi(x)$, this $\Phi$ can be interpreted by $\Phi f(x, y) = e'(x,y)f(\psi(x), \phi(x,y))$, with $\phi(x,\cdot)$ a measure preserving transformation for a.e. $x$. But since $T\mathbf 1 = \mathbf 1$, this $e'(x,y) = 1$ as well.  \\
	
	It remains to prove equimeasurability.  Let $\mathbf 1_{\mathbf B} = (\mathbf 1_{B_1},...,\mathbf 1_{B_n})$, and observe that since for a.e. $x$, $\phi(x,\cdot)$ is a measure preserving transformation inducing a lattice isometry over $L_q$, it follows that 
	$$N^q[\mathbf 1_{B_i}](x) =\mu(y: (x,y)\in B_i) = \mu(y:(x,\phi(x,y)) \in B_i),$$
	While
	\begin{align*}
	N^q[T\mathbf 1_{B_i}](x) =\mu(y: (\psi(x),\phi(x,y))\in B_i) \\
	 = \mu(y: (\psi(x),y)\in B_i) = N^q[\mathbf 1_{B_i}](\psi(x)).
	\end{align*}
	Thus for each $A = \prod_i A_i$ with $A_i\subseteq [0,1]$ measurable, since $\psi$ is also a measure preserving transformation, 
	\begin{align*}  \mu(x: (N^q[\mathbf 1_{\mathbf B}](x) \in A) = \mu(x: (N^q[\mathbf 1_{\mathbf B}](\psi (x)) \in A)  = \mu(x: (N^q[T\mathbf 1_{\mathbf B}](x) \in A),
	 \end{align*}
	and we are done.
	
\end{proof}

	 The following theorem describes a comparable equimeasurability property of certain copies of $L_pL_q$ in $L_p(L_q)$ for any $1\leq p \neq q <\infty$:

\begin{theorem}\label{t:blplq_equimeasurable}
	Let $1 \leq p \neq q < \infty$, and suppose that $f_i:E\rightarrow L_p(L_q)$ are lattice embeddings with $E\in B\K{p}{q}$ generated by a $(k,j)$-indexed collection of atoms $\mathbf e:= (e(k,j))_{k,j}$ with $1\leq k \leq n$ and $1\leq j \leq m_k$ as described in Proposition \ref{p:BLpLq-representation}.  Suppose also that $f(\sum_{k,j} e(k,j))= \mathbf 1 \cdot \|\sum e(k,j) \|.$  Then $(f_1(\mathbf e))$ and $(f_2(\mathbf e))$ are base-equimeasurable.  
\end{theorem}

\begin{proof}
	Let $\eta = \| \sum_{k,j} e(k,j) \|$, and note first that each $\frac{1}{\eta}f_i(e(k,j))$ is of the form $\mathbf 1_{A_i(k,j)}$ for some measurable $A_i(k,j) \subseteq [0,1]^2$.  Second,  $N^q[\mathbf 1_{A_i(k,j)}] (s) = \mu( A_i(k,j)(s) )$ with $A_i(k,j)(s) \subseteq [0,1]$ measurable for a.e. $s$, so by Proposition \ref{p:BLpLq-representation}, for each fixed $k$ and each $j, j'$, $\mu( A_i(k,j)(s) ) = \mu( A_i(k,j')(s) ) = \frac{1}{m_k}\mathbf 1_{A_i(k)} (s)$ with $A_i(1),...,A_i(n) \subseteq [0,1]$ almost disjoint. It follows that for each appropriate $k,j$, $\frac{1 }{\eta} = \frac{1}{m_k^{1/q}} \mu(A_i(k))^{1/p} $, so $\mu(A_i(k))  = \bigg(\frac{m_k^{1/q}}{\eta} \bigg)^p$.  \\
	
	To show equimeasurability, observe that for a.e. $t$, we have 
	$ N^q[\mathbf 1_{A_i(k,j)}](s) = \frac{1}{m_k} $ iff $s\in A_i(k)$, and 0 otherwise.  Let $\mathbf B \subseteq \prod_k \R^{m_k}$ be a measurable set. Note then that any $(k,j)$-indexed sequence $(N[f_i(\mathbf e)](s))$ is of the form $\mathbf {c^i_s} \in \prod_k \R^{m_k}$ with $c^i_s(k,j) = \bigg(\frac{1}{m_k}\bigg)^{1/q}$ for some unique $k$, and $c^i_s(k,j) = 0$ otherwise. It follows then that for some $I \subseteq {1,..., n}$,
	\[\mu(s: \mathbf{ c^i_s} \in \mathbf B) = \sum_{k\in I} \mu(A_i(k)) = \sum_{k\in I} \bigg(\frac{m_k^{1/q}}{\eta}\bigg). \]
	
	Since the above holds independent of our choice of $i$, we are done.
	
\end{proof}

\begin{remark}\label{r:BKpq_standard_representation}
	The above proof shows much more than base-equimeasurability for copies of $B\K{p}{q}$ lattices in $L_p(L_q)$.  Indeed, if $\mathbf 1\in E = < (e(k,j))_{k,j}> $ with $E \in B\K{p}{q}$, then each atom is in fact base-simple, and $\sum e(k,j) = \eta \cdot \mathbf 1$ where $\eta = (\sum_k m_k^{p/q} )^{1/p}$.  Furthermore, there exist measurable sets $A(1),...,A(n)$ partitioning $[0,1]$ with $\mu(A(k)) = \frac{m_k^{p/q}}{\eta^p}$ such that $N[e(k,j)] = \frac{\eta}{m_k^{1/q}} \mathbf 1_{A(k)}$. Based on this, we can come up with a "canonical" representation of $E$, with $e(k,j)\mapsto \eta\cdot \mathbf 1_{W_k \times V_{k,j}}$, where $$W_k = \big[\sum_{l=1}^{k-1} \mu(A(l)), \sum_{l=1}^k \mu(A(l)) \big] \text{, and } V_{k,j} =\bigg[\frac{j-1}{m_k}, \frac{j}{m_k} \bigg]. $$
	
	This canonical representation will become relevant in later results.
\end{remark}

Having characterized representations of lattice in $B\K{p}{q}$, we now move towards proving the AUH result.  Before the final proof, we use the following perturbation lemma. 

\begin{lemma}\label{l:approximate_with_full_support}
	Let $f:E\rightarrow L_p(L_q)$ be a lattice embedding of a lattice $E= <e_1,...,e_n>$.  Then for all $\varepsilon > 0$, there exists an embedding $g:E\rightarrow L_p(L_q)$ such that $g(E)$ fully supports $L_p(L_q)$ and $\|f-g\| < \varepsilon$.
\end{lemma}

\begin{proof}
	
	Let $M_k = supp\big(N[f(e_k)] \big) \backslash supp\big(N[f(\sum_1^{n-1}e_k)]\big)$. For each $e_k$, we will construct $e'_k$ disjoint from $f(E)$ with support in $M_k\times [0,1]$. Let $M'$ be the elements in $[0,1]^2$ disjoint from $f(E)$. Starting with $n = 1$, Observe that $M'$ can be partitioned by $M' \cap M_k \times [0,1]:= M'_k$. Let 
	
	$$\eta_k(x,y) = \varepsilon^{1/q} \frac{N[f(e_k)](x)}{\mu(M'_k(x))^{1/q}} \mathbf 1_{M'_k}(x,y). $$
	
	When $\mu(M'_k(x)) = 0$, let $\eta_k(x,y) = 0$ as well.  Now, let $g':E\rightarrow L_p(L_q)$ be the lattice homomorphism induced by $$g'(e_k) = (1-\varepsilon)^{1/q}f(e_k)\cdot \mathbf 1_{M_k} + \eta_n +f(e_k) \cdot \mathbf 1_{M_k^c}.$$
	First, we show that $g'$ is an embedding. Observe that for each $k$, 
	\begin{align*}
		N^q[g'(e_k)](x) = & \int \eta^q_k(x,y) + (1-\varepsilon)f(e_k)^q(x,y) \ dy \\
		= & \int \varepsilon \frac{N^q[f(e_k)](x)}{\mu(M'_k(x))} \cdot \mathbf 1_{M'_k}(x,y) + (1-\varepsilon)f(e_k)^q(x,y) \ dy \\
		= & \varepsilon N^q[f(e_k)](x) +(1-\varepsilon)\int f(e_k)^q(x,y) \ dy \\
		= & \varepsilon  N^q[f(e_k)](x) + (1-\varepsilon)  N^q[f(e_k)](x) =  N^q[f(e_k)](x).
	\end{align*}
	
	It easily follows that $g'(E)$ is in fact isometric to $f(E)$, and thus to $E$.  Furthermore, for every $k$, 
	
	\begin{align*}
		\|f(e_k) - g'(e_k)\| = &  \|\mathbf 1_{M_k} [ (1-(1-\varepsilon)^{1/q})f(e_k) + \eta_k ] \| \\
		\leq & (1-(1-\varepsilon)^{1/q}) + \varepsilon.
	\end{align*}
	The above can get arbitrarily small. \\
	
	Now, if $supp(N(\sum e_k)) = [0,1]$, let $g = g'$, and we are done.  Otherwise, let $\tilde M = \cup_k M_k$, and observe that $\sum g'(e_k)$ fully supports $L_p( \tilde M; L_q)$. Observe also that $L_p(L_q) = L_p(\tilde M; L_q) \oplus_p L_p(\tilde M^c; L_q)$.  However, both $L_p(\tilde M; L_q)$ and $L_p(\tilde M^c; L_q)$ are lattice isometric to $L_p(L_q)$ itself. So there exists an isometric copy of $E$ fully supporting $L_p(\tilde M^c; L_q)$. Let $e'_1,...,e'_n \in L_p(\tilde M^c; L_q)$ be the corresponding basic atoms of this copy, and let $g(e_i) = (1-\varepsilon^p)^{1/p}g'(e_i) + \varepsilon\cdot e'_n$. Then for $x\in E$, 
	
	$$\|g(x)\|^p = (1-\varepsilon) \|g'(x) \|^p+ \varepsilon \| x \|^p = \|x\|^p. $$
	
	Using similar reasoning as in the definition of $g'$, one also gets $\|g - g'\| < (1-(1-\varepsilon)^{1/p}) + \varepsilon$, so $g$ can also arbitrarily approximate $f$.  
	
\end{proof}

Observe that the lemma above allows us to reduce the approximate homogeneity question down to cases where the copies of a $B\K{p}{q}$ lattice fully support $L_p(L_q)$.  Combined with Proposition \ref{p:homogeneity-one-element}, we can further reduce the possible scenarios to cases where for each $i$, $f_i(x) = \mathbf 1$ for some $x\in E$.  It turns out these reductions are sufficient for constructing a lattice automorphism that makes the homogeneity diagram commute as desired:

\begin{theorem}\label{t:blplq_full_homogeneity}
	Suppose $1\leq p\neq q < \infty$, and for $i = 1,2$, let $f_i:E \rightarrow L_p(L_q)$ be a lattice embedding with $E := < (e(k,j) )_{k,j}> \in B\K{p}{q}$ and $1\leq k \leq n$ and $1\leq j \leq m_k$. Suppose also that each $f_i(E)$ fully supports $L_p(L_q)$.  Then there exists a lattice automorphism $\phi$ over $L_p(L_q)$ such that $\phi \circ f_1 = f_2$.  
\end{theorem}	

\begin{proof}
	
	Let $\eta = \| \sum_{k,j} e(k,j) \|$; by Proposition \ref{p:homogeneity-one-element}, we can assume that for both $i$'s, we have $f_i(\sum_{k,j} e(k,j)  ) = \eta \cdot \mathbf 1$. For notation's sake, let $e_i(k,j):=f_i(e(k,j)) $. By Proposition \ref{p:BLpLq-representation}, for each $i$ there exist mutually disjoint sets $A_i(1),...,A_i(n)$ partitioning $[0,1]$ such that for each $1\leq j \leq m_k$, $supp(N[e_i(k,j)]) = A_i(k)$. In addition, for the sets $A_i(k,1),...,A_i(k,m_k)$, where $A_i(k,j):= supp(e_i(k,j))$, partition $A_i(k) \times [0,1]$. It follows also from the statements in Remark \ref{r:BKpq_standard_representation} that $\mu(A_1(k)) = \mu(A_2(k))$ for each $k$ and $N^q[e_i(k,j)](x) = \frac{\eta^q}{m_k} \mathbf 1_{A_i(k)}(x)$. \\

	To prove the theorem, it is enough to generate lattice automorphisms $\phi^i$ mapping each band $\bob(e_i(k, j))$ to a corresponding band $\bob(\mathbf 1_{W_k \times V_{k,j}})$ where $W_k$ and $V_{k,j}$ are defined as in Remark \ref{r:BKpq_standard_representation}, with $\mathbf 1_{A_i(k,j)} \mapsto \mathbf 1_{W_k \times V_{k,j}}$.\\
	
	 To this end, we make a modified version of the argument in \cite[Proposition 2.6]{henson11} and adopt the notation in Proposition \ref{p:homogeneity-one-element}: construct lattice isometries $\psi^i_{k,j}$ from $L_p(A_i(k)); L_q(V_{k,j}))$ to $\bob(e^i_{k,j})$ with 
	\[ \psi^i_{k,j}(f)(x,y) = f\bigg(x, \big( \widetilde {\mathbf 1}_{A_i(k,j)}\big)_x(y)_q + \frac{j-1}{m_k} \bigg) \mathbf 1_{A_i(k,j)} (x,y)\]
By similar reasoning as in the proof of Proposition \ref{p:homogeneity-one-element}, $\psi^i_{k,j}$ is a lattice embedding.  Surjectivity follows as well. Indeed, since $N^q[\mathbf 1_{A_i(k,j)} ](x) = \frac{1}{m_k}$, for a.e. $x\in A_i(k)$ the function $ \big( \widetilde {\mathbf 1}_{A_i(k,j)} \big)_x(y)_q + \frac{j-1}{m_k}$ matches $[0,1]$ continuously to $V_{k,j}$ with $supp(e_i(k,j)(x,\cdot))$ mapped a.e. surjectively to $V_{k,j}$. So $\psi^i_{k,j}$'s image is dense in $\bob (e_i(k,j))$. \\
	
Observe that $\psi^i_{k,j}$ also preserves the random norm $N$ along the base (that is: $N[f] = N[\psi^i_{k,j}(f)]$. Resultantly, the function $\psi^i_k := \oplus_j \psi^i_{j,k}$ mapping $L_p(A_i(k), L_q(0,1))$ to $\oplus_j \bob(e_i(k,j))$ is also a lattice automorphism.  Indeed, for $f = \sum_1^{m_k} f_j$ with $f_j \in \bob(e_i(k,j))$, one gets

\begin{align*}
	\| \psi^i_k(f) \| & = \bigg|\bigg| N[\sum_j \psi_{k,j}^i(f_j)] \bigg|\bigg|_p = \bigg|\bigg| \big( \sum_j N^q[ \psi^i_{k,j}(f_j)] \big)^{1/q} \bigg|\bigg|_p \\
	& = \bigg|\bigg| \big( \sum_j N^q[ f_j] \big)^{1/q} \bigg|\bigg|_p = \bigg|\bigg| N[\sum_j f_j] \bigg|\bigg|_p=  \| f \|
\end{align*}
Now let $\psi^i = \oplus_k \psi^i_k$, and observe that given $f = \sum_1^n f_k$ with $f_k \in L_p(A_i(k), L_q(0,1))$, since the $f_k$'s are base disjoint, we have 

$$\|\psi^i f\|^p = \sum_1^n \| \psi^i_k f_k \|^p =  \sum_1^n \| f_k \|^p = \| f\|^p. $$

  Thus $\psi^i$ is a lattice automorphism over $L_p(L_q)$ mapping  each $1_{A_i(k) \times V_{k,j}}$ to $\mathbf 1_{A_i(k,j)}$.  \\
  
  Use \cite[Lemma 3.3]{fer20} to construct a lattice isometry $\rho_i: L_p \rightarrow L_p$ such that for each $k$, $\rho_i (\mathbf 1_{W_k}) = \mathbf 1_{A_i(k)}$. By \cite[Ch. 11 Theorem 5.1]{banach1932theorie} this isometry is induced by a measure preserving transformation $\bar{\rho_i}$ from [0,1] to itself such that $\rho^i(f)(x) = f(\bar{\rho_i}(x))$. It is easy to show that $\rho_i$ induces a lattice isometry with $f(x,y) \mapsto f(\bar{\rho_i}(x), y)$.  In particular, we have $N[\rho_i f](x) = N[f](\bar{\rho_i}(x))$ , and  $\rho_i(\mathbf 1_{W_k\times V_{k,j}}) = \mathbf 1_{A_i(k)\times V_{k,j}}$, now let $\phi^i(f) = (\psi^i\circ \rho^i)(f)$, and we are done.
	
\end{proof}

Using the above, we can now show:

\begin{theorem}\label{t:blplq_finite_homogeneity}
	For $1\leq p\neq q <\infty$, the lattice $L_p(L_q)$ is AUH for the class $B\K{p}{q}$. 
\end{theorem}

\begin{proof}
	Let $f_i:E\rightarrow L_p(L_q)$ as required, and suppose $\varepsilon > 0$.  use Lemma \ref{l:approximate_with_full_support} to get copies $E'_i$ of $f_i(E)$ fully supporting $L_p(L_q)$ such that for each atom $e_k \in E$ and corresponding atoms $e^i_k \in E'_i$, we have $\|f_i(e_k) - e^i_k \| < \varepsilon/2$. now use Theorem \ref{t:blplq_full_homogeneity} to generate a lattice automorphism $\phi$ from $L_p(L_q)$ to itself such that $\phi(e^1_k) = e^2_k$.  Then 
	$$\| \phi(f_1(e_k)) - f_2(e_k)) \| \leq \|\phi(f_1(e_k) - e^1_k ) \| +\| e^2_k -f_2(e_k) \|< \varepsilon $$
\end{proof}

\begin{remark}\label{r:AUH_must_be_atomless}
	Observe that the doubly atomless $L_p(L_q)$ space is unique among separable $BL_pL_q$ spaces that are AUH over $B\K{p}{q}$. Indeed, this follows from the fact that such a space must be doubly atomless to begin with: let $E$ be a one dimensional space generated by atom $e$ and suppose $X$ is not doubly atomless.  Suppose also that $E$ is embedded by some $f_1$ into a part of $X$ supported by some $L_p$ or $L_q$ band, and on the other hand is embedded by some $f_2$ into $F:= \ell_p^2(\ell_q^2)$ with $f_2(e)$ a unit in $F$. 	Then one cannot almost extend $f_1$ to some lattice embedding $g:F\rightarrow X$ with almost  commutativity.
\end{remark}

One can also expand this approximate ultrahomogeneity to separable sublattices with a weaker condition of almost commutativity in the diagram for generating elements: for any $BL_pL_q$ sublattice $E$ generated by elements $<e_1,...,e_n>_L$, for any $\varepsilon > 0$, and for all lattice embedding pairs $f_i:E\rightarrow L_p(L_q)$, there exists a lattice automorphism $g:L_p(L_q) \rightarrow L_p(L_q)$ such that for all $j = 1,...,n$, $\|g(f_2(e_j)) - f_1(e_j)\| <\varepsilon$. \\

\begin{theorem}\label{t:finite-generated-blplq-homogeneity}
	For all $1\leq p\neq q<\infty$, The lattice $L_p(L_q) $ is AUH for the class of finitely generated $BL_pL_q$ lattices.
\end{theorem}
\begin{proof}
	Let $E = <e_1,...e_n>_L$, and let $f_i:E\rightarrow L_p(L_q)$ be lattice embeddings. We can assume that $\|e_k\| \leq 1$ for each $1\leq i\leq n$.  By Proposition \ref{p:finite_representation2}, $E$ is the inductive limit of lattices in $B\K{p}{q}$. Given $\varepsilon > 0$, pick a $B\K{p}{q}$ lattice $E' = <e'_1,...,e'_m> \subseteq E$ such that for each $e_k$, there is some $x_k \in B(E')$ such that $\| x_k - e_k \| < \frac{\varepsilon}{3}$.  Each $f_i|_{E'}$ is an embedding into $L_p(L_q)$, so pick an automorphism $\phi$ over $L_p(L_q)$ such that $\| \phi \circ f_1|_{E'} - f_2|_{E'} \| < \frac{\varepsilon}{3}$.  Then
	$$\| \phi f_1 (e_k) -f_2(e_k)  \| \leq \|\phi f_1(e_k - x_k) \| + \| \phi f_1(x_k) -f_2(x_k)\| + \| f_2(x_k - e_k) \| < \varepsilon.$$
\end{proof}

We can also expand homogeneity to include not just lattice embeddings but also disjointness preserving linear isometries, that is, if embeddings $f_i:E\rightarrow L_p(L_q)$ are not necessarily lattice homomorphisms but preserve disjointness, then there exists a disjointness preserving linear automorphism $\phi$ over $L_p(L_q)$ satisfying almost commutativity:

\begin{corollary}\label{c:homogeneity_disjointness}
	$L_p(L_q)$ is AUH over finitely generated sublattices in $BL_p(L_q)$ with disjointness preserving embeddings.
\end{corollary}

\begin{proof}
	Use the argument in \cite[Proposition 3.2]{fer20} to show that $L_p(L_q)$ is disjointness preserving AUH over $B\K{p}{q}$.  From there, proceed as in the argument in Theorem \ref{t:finite-generated-blplq-homogeneity} to extend homogeneity over $B\K{p}{q}$ to that over $BL_pL_q$.   
\end{proof}

\section{Approximate Ultrahomogeneity of $L_p(L_q)$ when $p/q\notin \N$}\label{s:equimeasurability}
The above results largely focused approximate ultrahomogeneity over $BL_pL_q$ lattices. What can be said, however, of \textit{sublattices} of $L_pL_q$ spaces?   The answer to this question is split into two cases: first, the cases where $p/q \notin \N$, and the second is when $p/q \in \N$.  We address the first case in this section. It turns out that if $p/q \notin \N$, then $L_p(L_q)$ is AUH for the class of its finitely generated sublattices. The argument involves certain equimeasurability properties of copies of fixed finite dimensional lattices in $L_p(L_q)$. Throughout, we will refer to the class of sublattices of spaces in $B\K{p}{q}$ as simply $\K{p}{q}$, and let $\overline{\K{p}{q}}$ be the class of finitely generated sublattices of $L_p(L_q)$. \\

The following result appeared as \cite[Proposition 3.2]{henson11}, which is a multi-dimensional version based on Raynaud's proof for the case of $n = 1$ (see \cite[lemma 18]{raynaud1986sous}). The approach taken here is a multi-dimensional version of the proof of Lemma 2 in \cite{Linde1986UniquenessTF}. 

\begin{theorem}\label{theorem:lindeExtension}
	Let $r = p/q \notin \N$, and suppose $f_i:E \rightarrow L_p(L_q) $ are lattice isometric embeddings with $E = <e_1,...,e_n>$. Suppose also that $f_1(x) = f_2(x) = \mathbf 1$ for some $x\in E_+$. Then $f_1(\mathbf e)$ and $f_2(\mathbf e)$ are base-equimeasurable.
\end{theorem}

 Throughout the proof, let $\mu$ be a measure in some interval $I^n \subseteq C:=\R^n_+$. To this end, we first show the following:

\begin{lemma}\label{lemma:LindeExtention}
	Suppose $0 < r \notin \N$, and $\alpha, \beta$ are positive finite Borel measures on $C$ such that for all $\mathbf v \in C$ with $v_0 > 0$,  
	\[ \int_{C} |v_0 + \mathbf v \cdot \mathbf z |^r \ d \alpha(\mathbf z) = \int_{C} |v_0 + \mathbf v \cdot \mathbf z |^r \ d \beta(\mathbf z) <\infty. \]
	
	Then $\alpha = \beta$.
\end{lemma}

\begin{proof}
	
	It is equivalent to prove that the signed measure $\nu:=\alpha - \beta = 0$. 	First, observe that since $|\nu| \leq \alpha +\beta$, and for any $\mathbf v \geq 0$, $\int |v_0 + \mathbf v\cdot \mathbf z|^r \ d|\nu|(\mathbf z) < \infty$. \\
	
	Now, we show by induction on polynomial degree that for all $ k \in \N$, $\mathbf v \geq 0$, and for all multivariate polynomials $P(\mathbf z)$ of degree $k'\leq k$,
	
	\[ * \int_{C} |v_0 + \mathbf v \cdot \mathbf z |^{r-k} P(\mathbf z) \ d \nu(\mathbf z) = 0. \]
	
	This is true for the base case $k = 0$ by assumption.  Now assume it is true for $k \in \N$ and let $k':= \sum l_i \leq k$ with $\mathbf l \in \N^n$.  For notational ease, let $\mathbf z^{\mathbf l} = z_1^{l_1}...z_n^{l_n}$.  Then for each $v_i$ and $0 < t < 1$, 
	
	\[ \int_{R^n_+} \mathbf z^{\mathbf l} \frac{(v_0 + \mathbf v\cdot \mathbf z +  z_i t)^{r-k} - (v_0 + \mathbf v\cdot \mathbf z)^{r-k} }{t} \ d \nu(\mathbf z) = 0. \]
	
	Now, if $k+1< r$ and $t\in (0,1)$, then
	
	\begin{align*}  & \bigg| \mathbf z^{\mathbf l}\frac{(v_0 + \mathbf v\cdot \mathbf z + z_i t)^{r-k} - (v_0 + \mathbf v\cdot \mathbf z)^{r-k} }{t} \bigg|
		\leq \mathbf z^{\mathbf l} z_i (r- k) (v_0 +\mathbf v\cdot \mathbf z + v_i)^{r-k-1} \\
		\leq & \frac{r-k}{\mathbf v^{\mathbf l} v_i}(v_0 + \mathbf v \cdot \mathbf z +v_i)^r
	\end{align*}
	
	Since in this case, $0 < r - k -1 < r $ and $|\nu| < \infty$, the right hand side must also be $|\nu|$-integrable. On the other hand,  If $k+1 > r$,  then we have 
	
	\[ \bigg|  \mathbf z^{\mathbf l}\frac{(v_0 + \mathbf v\cdot \mathbf z + v_i t)^{r-k} - (v_0 + \mathbf v\cdot \mathbf z)^{r-k} }{t} \bigg|  < |r - k|\frac{ v_0^{r }}{\mathbf v^{\mathbf l} v_i} \]
	which is also $|\nu|$-integrable. So now we apply Lebesgue's differentiation theorem over $v_i$ to get, for any $k \in \N$ and for each $1 \leq i \leq n$:
	
	\[  \int_{C} \mathbf z^{\mathbf l} z_i|v_0 + \mathbf v\cdot \mathbf z|^{r-k-1} \ d\nu(\mathbf z) = 0,  \]
	
	since $r\notin \N$. A similar argument, deriving over $v_0$, can be made to show that 
	
	\[ \int_{C} |v_0 + \mathbf v\cdot \mathbf z|^{r-k-1} \ d\nu(\mathbf z) = 0 \]
	
	One can make linear combinations of the above, which implies line $*$. \\
	
	Now for fixed $\mathbf v > 0$, $v_0 >0$ we define a measure $\Lambda$ on $C$, where for measurable $B \subseteq \R^n_+,$
	
	\[ \Lambda( B ) = \int_{\phi^{-1}(B)} |v_0 + \mathbf v \cdot \mathbf z|^r \ d\nu(\mathbf z).  \]
	
	where $\phi(\mathbf z)= \frac{1}{v_0 +\mathbf v \cdot \mathbf z} \mathbf z $. It is sufficient to show that $\Lambda = 0$.  Observe first that $\phi$ is continuous and injective; indeed, if $\phi(\mathbf z) = \phi(\mathbf w)$, then it can be shown that $\mathbf v\cdot \mathbf w = \mathbf v\cdot \mathbf z$. Thus  $\frac{\mathbf w}{v_0 +\mathbf v \cdot \mathbf z} = \frac{\mathbf z}{v_0 +\mathbf v \cdot \mathbf z}$, implying that $\mathbf w = \mathbf z$.  Resultantly,  $\phi(B)$ for any Borel $B$ is also Borel, hence we will have shown that for any such $B$, $\nu(B) = 0$ as well, so $\nu = 0$.\\
	
	Observe that by choice of $\mathbf v >0$ and and since $(v_0 +\mathbf v\cdot \mathbf z) > 0$ for all $\mathbf z\in \R^n_+$, have
	
	\[ |\Lambda| (B) = \int_{\phi^{-1}(B)} |v_0 +\mathbf v\cdot \mathbf z| ^r \ d|\nu|(\mathbf z). \]
	
	Using simple functions and the definition of $\Lambda$, one can show both that for each $i$, we have 
	\begin{align*} ** \ \ m_i(k):= \int_C w_i^k \ d |\Lambda|(\mathbf w) = \int_C (v_0 +\mathbf v\cdot \mathbf z)^{r-k}z_i^k \ d\|\nu|(\mathbf z) <\infty 
	\end{align*}	 
	and also that 
	
	\[  \int_C w_i^k \ d \Lambda(\mathbf w) = \int_C |v_0 +\mathbf v\cdot \mathbf z|^{r-k}z_i^k \ d \nu(\mathbf z) = 0,   \]
	
	More generally, if $k = \sum_i l_i$, then 
	\begin{align*} \int_C \mathbf w^{\mathbf l} \ d\Lambda(\mathbf w) = \int_C \mathbf z^{\mathbf l} |v_0 +\mathbf v \cdot \mathbf z|^{r-k} \ d\nu(\mathbf z) = 0,
	\end{align*}
	So it follows that $ \int_C P(\mathbf w) \ d\Lambda(\mathbf w) = 0$ for all polynomials $P(\mathbf w)$.\\
	
	Now if $k > r$ and  $\nu \neq 0$, since $v_i > 0$, we then we have
	\begin{align*} 
		m_i(k) & =  \int_C |v_0 +\mathbf v\cdot \mathbf z|^{r-k}z_i^{k} \ d\|\nu|(\mathbf z) \\ 
		&\leq  \int_C |v_0+\mathbf v\cdot \mathbf z|^r v_i^{-k} \ d|\nu|(\mathbf z) \leq v_i^{-k}|\Lambda|(C)  <\infty \end{align*}
	
	so 
	\begin{align*}
		m_i(k)^{-1/2k} \geq v_i^{1/2} |\Lambda|(C)^{-1/2k}
	\end{align*}
	
	Thus for each $1 \leq i \leq n$, $ \sum_k m_i(k)^{-1/2k} = 0$.  So by \cite[Theorem 5.2]{jeu2001}, $|\Lambda|$ is the unique positive measure over $C$ with moment values $m_i(k)$.  Since $|\Lambda| +\Lambda$ yields the same values, and by **, $\int_C P(\mathbf w) \ d (|\Lambda| + \Lambda)(\mathbf w) = \int_C P(\mathbf w) \ d|\Lambda|(\mathbf w)$,  it follows that $\Lambda = 0$, so $\nu = 0$.

\end{proof}

Now we are ready to prove Theorem \ref{theorem:lindeExtension}.

\begin{proof}
	
	For simplicity of notation, let $F^i_j = N^q[f_i(e_j)]$ and $I = [0,1]$. By definition of $N$, the support of $F^i_j$ as well as of $\mu$ is the unit interval. Define positive measures $\alpha_j$ by 
	\[ \alpha_i(B) = \mu(\{ t \in I : \mathbf F^i(t) \in B \}) = \mu((\mathbf F^{i})^{-1}(B)). \]
	
	Now, for any measurable $B  \subseteq C$, we have 
	\begin{align*}
		\int_C \mathbf 1_B(\mathbf z) \ da_i(\mathbf z) = \alpha_i(B) = \mu((\mathbf F^i)^{-1}(B)) =\int_I (\mathbf 1_B\circ\mathbf  F^i) (t) \ dt
	\end{align*}
	so for any simple function $\sigma$ over $C$,
	\[ \int_C \sigma(\mathbf z) \ d \alpha_i = \int_0^1 \sigma \circ \mathbf F^i (t) \ dt \]
	
	Using simple functions to approximate $|v_0+\mathbf v\cdot \mathbf z|^r$, and given that $|v_0+\mathbf v\cdot \mathbf z|^r$ is in $L_1(C,\mu)$ and the support of $\mu$ is the unit interval, it follows that
	
	\[ \int_C |1+\mathbf v\cdot \mathbf z|^r \ d \alpha_i(\mathbf z) = \int_0^1 |1+\mathbf v \cdot \mathbf F^i(t)|^r \ dt.  \]
	
	It is sufficient now to show that for all $\mathbf v \in \R^n_+$,
	
	\[ \int_0^1 |1+\mathbf v \cdot \mathbf  F^1(t)|^r \ dt = \int_0^1 |1+\mathbf v \cdot \mathbf F^2(t)|^r \ dt.  \].
	
	For $i,j$ and $s\in [0,1]$, let $M_i^j = \{(s,t): (s,t) \in supp(f_i(e_j))\}$, and let $M_i^j(s) = \{t: (s,t) \in M_i^j \}$. By assumption, $x =\sum_j x_j e_j$ with $x_j > 0$, so  $\mathbf 1 = N^q[f_i(x)] = \sum_j x_j^q F^i_j$.   Therefore, since each $f_i$ is an embedding, for all $\mathbf c \in \R^n_+$, 
	
	\begin{align*}
		\|\sum_j c_j e_j  \|^p = & \bigg|\bigg| \big(\sum_j c_j^q F^i_j(s) \big)^{1/q} \bigg|\bigg|_p \\
		= & \bigg|\bigg| \big( \mathbf 1 + \sum_j (c_j^q -x_j^q) F^i_j(s) \big)^{1/q} \bigg|\bigg|_p \label{line:homogeneity-with-equimeasurability}
	\end{align*} 
	Let $v_j:= c^q_j - x_j^q$: then in particular it follows that for all $\mathbf v \geq  0$, we have
	
	\[ \int_0^1  \bigg( 1 +  \mathbf v \cdot \mathbf F^1(s) \bigg)^{p/q} \ ds  = \int_0^1  \bigg( 1 +  \mathbf v \cdot \mathbf F^2(s) \bigg)^{p/q} \ ds.\]
	
	By Lemma \ref{lemma:LindeExtention}, we can conclude that $\alpha_1 = \alpha_2$, so $\mathbf F^1$ and $\mathbf F^2$ are equimeasurable. 
\end{proof}

Using Theorem \ref{theorem:lindeExtension}, we can uniquely characterize lattices in $\K{p}{q}$ in a way that parallels Proposition \ref{p:BLpLq-representation}.

\begin{theorem}\label{t:extend_to_Blplq}
	Suppose that $p/q \notin \N$, and let $E\subseteq L_p(L_q)$ with $E = <e_1,...,e_m>$.  Then the following hold:
	
	\begin{itemize}
		\item $E \in \K{p}{q}$ iff there exist mutually disjoint measurable functions $\phi(k,j) \in S(L_p(L_q))_+$, with $1\leq j \leq n$ and $1\leq k \leq L$ such that for each $j$, $e_j \in < (\phi(k,j))_k > = \ell_p^n$, and $<(\phi(k,j))_{k,j}> \in B\K{p}{q}$.
		
		\item Suppose $f_i:E\rightarrow L_p(L_q)$ is a lattice embedding with $i = 1,2$ and $E\in \K{p}{q}.$  Then there exist embeddings $f'_i: E' \rightarrow L_p(L_q)$ extending $f_i$ such that $E' \in B\K{p}{q}$.
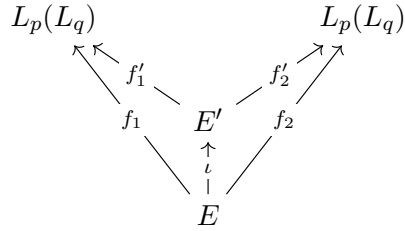
		\begin{center}
			\begin{tikzcd}
				L_p(L_q) & \quad & L_p(L_q) \\
				\quad & E' \arrow{lu}[description]{f'_1}  \arrow{ru}[description]{f'_2 } & \quad\\
				\quad & E  \arrow{luu}[description]{f_1}  \arrow{ruu}[description]{f_2}  \arrow{u}[description]{\iota}  & \quad
			\end{tikzcd}  
		\end{center}	
		
	\end{itemize}
\end{theorem}

\begin{proof}
	
	For part 1, clearly the reverse direction is true.  To prove the main direction, we can suppose that $E$ fully supports $L_p(L_q)$.  If not, recall that the band generated by $E$ is itself doubly atomless, and hence is lattice isometric to $L_p(L_q)$ itself. Thus, if under these conditions, there is a $BL_pL_q$ sublattices extending $E$ as in the statement of the theorem, it will also be the case in general. \\
	
	By Proposition \ref{p:homogeneity-one-element}, we can also suppose that $\sum_j e_j = \eta \cdot \mathbf 1$.  Now by assumption, since $E \in \K{p}{q}$, then there is an embedding $\psi:E\rightarrow \widetilde E \in B\K{p}{q}$ such that each $\psi(e_j) = \sum_k x(k,j)\tilde e(k,j)$, with $1\leq k \leq m'_k$. Without loss of generality we may also drop any $\tilde e(k,j)$'s disjoint from $\psi(E)$ and assume that $\psi(E)$ fully supports $\widetilde E$. Now $\widetilde E$ is a $B\K{p}{q}$ lattice admitting a canonical representation in $L_p(L_q)$ as described in Theorem \ref{t:blplq_equimeasurable} and Remark \ref{r:BKpq_standard_representation}.\\
	
	So we can assume that $\psi$ embeds $E$ into $L_p(L_q)$ in such a way that $\psi(E)$ fully supports it and each $\psi(e_j)$ is both simple and base-simple. Now, use Proposition \ref{p:homogeneity-one-element} to adjust $\psi$ into an automorphism over $L_p(L_q)$ such that $\psi(\sum e_j) = \eta\cdot  \mathbf 1$ in a way that preserves both simplicity and base-simplicity. By Theorem \ref{theorem:lindeExtension}, $\psi(\mathbf e)$ and $\mathbf e$ are base-equimeasurable.  Since the $\psi(k,j)'s$ are base-simple, there exist tuples $\mathbf {s^1},...,\mathbf{ s^L} \in \R^m$ such that for a.e. $t\in [0,1]$, there is some $k \leq L$ such that $N[\mathbf e](t) = \mathbf {s^k}$. By equimeasurability, the same is true for $N[{\psi(\mathbf e)}](t)$.  \\
	
	Let $\mathbf {S^k} = \{ t: N[\mathbf e](t) = \mathbf {s^k } \}$, and let $S^k_j = \mathbf {S^k }\times [0,1] \cap supp(e_j) $. Let $\mathbf{\overline S^k} =  \{ t: N[\psi(\mathbf e)](t) = \mathbf {s^k } \}$ with $\overline S^k_j$ defined similarly. Note that each $\mathbf 1_{S^k_j}$  is also base-characteristic, as $N[\mathbf 1_{S^k_j}] = c^k_j \mathbf 1_{\mathbf {S^k}}$ for some $c^k_j > 0$, so for fixed $k$ and for any $j,j' \leq m_k$, we must have that $N[\mathbf 1_{S^k_j}]$ and $N[\mathbf 1_{S^k_{j'}}]$ are scalar multiples of each other. Thus for each appropriate pair $(k,j)$ with $s^k_j > 0$, define $\phi(k,j)$ by $\frac{\mathbf 1_{S^k_j}}{\|\mathbf 1_{S^k_j}\|}$.  By definition of $\mathbf{ S^k}$, for any $k \neq k'$ and any appropriate $j,j'$, $\phi(k,j)$ and $\phi(k',j')$ are fiber-disjoint, and $N[\phi(k,j)] = N[\phi(k,j')]$.  Thus by Proposition \ref{p:BLpLq-representation}, $ <(\phi(k,j))_{k,j}> \in B\K{p}{q}$. \\
	
	To prove part 2, Observe first that we have already essentially proven part 2 in the case that $f_1 =Id$ and $f_2 = \psi$.  To show the general case, we first assume that for each $i$, $\sum f_i(e_j)$ maps to $\mathbf 1$. Now, by Theorem \ref{theorem:lindeExtension}, $f_1(\mathbf e)$ and $f_2(\mathbf e)$ are also base-equimeasurable, but by the procedure for part 1, we also know that each $f_i(e_j)$ is also base-simple. Define $\mathbf{s^1},...,\mathbf {s^L}$ as above, and Let $ \mathbf {S^k}(i) = \{ t: N[f_i(\mathbf e)](t) = \mathbf{s^k}  \}$. Define similarly $S^k_j(i)$ and the associated characteristic functions $\phi_i(k,j)$ for appropriate pairs $k,j$ such that $1\leq k \leq l$ and $s^k_j: = \|\phi_i(k,j) \wedge f_i(e_j)\| > 0$.  Note first that $$f_i(e_j) = \sum_{k: s_k(j) >0} s^k_j \phi_i(k,j).$$ Second, observe that by equimeasurability, the eligible pairs $(k,j)$ are the same for $i=1,2$.  Let $E'_i = <(\phi_i(k,j))_{k,j} >$. Clearly $E'_i \in B\K{p}{q}$, and since the eligible pairs $(k,j)$ are the same, $E'_1$ and $E'_2$ are isometric to each other.  Let $E'$ be one of the $E'_i$'s and let $f'_i:E' \rightarrow L_p(L_q)$ be the expected embedding mapping $E'$ to $E'_i$, and we are done.
\end{proof}

From here, we can now easily extend Theorem \ref{t:blplq_full_homogeneity} to lattices in $\K{p}{q}$:

\begin{corollary}\label{l:homogeneity-with-equimeasurability}
	Suppose $p/q \notin \N$ and suppose $f_i:E\rightarrow L_p(L_q)$ are lattice embeddings from $E\in K_{p,q}$ with $f_i(E)$ fully supporting $L_p(L_q)$.  Then there exists a lattice automorphism $\phi$ over $L_p(L_q)$ such that $f_2 = \phi \circ f_1$. 
\end{corollary}  

\begin{proof}
	
	Use Theorem \ref{t:extend_to_Blplq} to generate a $B\K{p}{q}$ lattice $E'$ containing $E$ and lattice embeddings $f'_i:E' \rightarrow L_p(L_q)$ such that $f'_i |_E = f_i$.  Clearly each $f'_i(E')$ fully supports $L_p(L_q)$.  Now apply Theorem \ref{t:blplq_full_homogeneity} to generate an automorphism $\phi$ over $L_p(L_q)$ with $\phi \circ f'_1 = f'_2$.  Clearly $\phi \circ f_1 = f_2$ as well.
	
\end{proof}

When $p/q \notin \N$, using Theorem \ref{t:extend_to_Blplq}, we can show that the same holds with the more general class  $\K{p}{q}$.  However, we can make an even stronger claim by showing that homogeneity holds for any finite dimensional sublattice of $L_p(L_q)$.  This is done using the following result, which gives a standard way of approximating finite dimensional sublattices of $L_p(L_q)$ with lattices in $\K{p}{q}$.

\begin{lemma}\label{l:approximate_with_simple_space}
	Suppose $p/q \notin \N$, and let $f_i:E\rightarrow L_p(L_q)$ be embeddings with $E = <e_1,...,e_n>$.  Then for all $\varepsilon > 0$, there exists a $\K{p}{q}$ lattice $E' = <e'_1,...,e'_n>$ and embeddings $g_i:E'\rightarrow L_p(L_q)$ such $g_i(E')$ fully supports $L_p(L_q)$ and for each $n$, $\|f_i(e_n) - g_i(e'_n) \| < \varepsilon$.
\end{lemma}

\begin{proof}
	We can assume each $f_i(E)$ fully supports $L_p(L_q)$: given $\varepsilon >0$, use Lemma \ref{l:approximate_with_full_support} to get copies of $E$ sufficiently close to each $f_i(E)$ with full support.  We then also assume that $f_i(\sum_1^n e_k) = \mathbf 1$ using Proposition \ref{p:homogeneity-one-element}.\\
	
	By Theorem \ref{theorem:lindeExtension}, $f_1(\mathbf e)$ and $f_2(\mathbf e)$ are base-equimeasurable.  In particular, given any measurable $C \in \R^n$, one has $\mu(t:N[f_1(\mathbf e)](t) \in C) = \mu(t: N[f_2(\mathbf e)](t) \in C)$.  Now pick an almost disjoint partition $C_1,...,C_m$ of $S(\ell_1^n)$, where each $C_i$ is closed, has relatively non-empty interior, and is of diameter less than $\frac{\varepsilon}{2n}$. Let $D_k^i = \{ t:N[f_i(\mathbf e)](t) \in C_i\backslash \cup_j^{i-1}C_j\}$.  Then by equimeasurability, $\mu (D_k^1) = \mu(D_k^2)$.   For each $k$, pick some $\mathbf {s^k} = (s^k_1,...,s^k_n) \in C_k$, and for each $x\in D_k^i$, let 
	$$\overline{e}^i_j(x,y) = \frac{s^k_j}{N[f_i(e_j)](x)}f_i(e_j)(x,y).$$
	
	Observe that $\|\sum_j\overline{e}^i_j - \sum_jf_i(e_j) \| < \varepsilon$, and $N[\overline{e}^i_j](x) = s_j^k$ for $x\in D^i_k$.  
	
	Consider now the lattice $E'= <\overline e^1_j,..., \overline e^1_n> $. Now, for any linear combination $\sum a_j \overline e^i_j$, we have, as in the argument in Proposition \ref{p:homogeneity-one-element}, that 
	
	\begin{align*}
		\| \sum a_j \overline e^i_j \|^p = \sum_k^M (\sum_j (a_j s^k_j)^q)^{p/q}
	\end{align*}
	
	implying that $\| \sum a_j \overline e^1_j \| = \| \sum a_j \overline e^2_j \|$. It follows both that $E'$ embeds into $\ell_p^M(\ell_q^n)$, implying that it  is a $\K{p}{q}$ lattice, and it is isometric to the lattice generated by the $\overline e^2_j$'s.  Let $e'_j = \overline e^1_j$, and define $g_i:E' \rightarrow L_p(L_q)$ as the maps generated by $g_i(e'_j) = \overline e^i_j$.  Clearly these are lattice embeddings and $\|f_i(e_j) - g_i(e'_j)\| < \varepsilon$.
	
\end{proof}

\begin{theorem}\label{t:homogeneous-finite-dim}
	For all $1\leq p, q < \infty$ with $p/q \notin \N$, the lattice $L_p( L_q)$ is AUH for the class of finite dimensional sublattices of $L_pL_q$ lattices.
\end{theorem}

\begin{proof}
	It is sufficient to show that the result is true over generation by basic atoms.  Let $f_i:E \rightarrow L_p(L_q)$ be two embeddings with $E = <e_1,...,e_n>$. Use Lemma \ref{l:approximate_with_simple_space} to find $g_i:E' \rightarrow L_p(L_q)$, with $E' := <e'_1,...,e'_n> \in \K{p}{q}$,  $\| g_i(e'_k) - f_i(e_k) \| < \varepsilon/2$,  and each $g_i(E')$ fully supporting $L_p(L_q)$.  Then by Lemma \ref{l:homogeneity-with-equimeasurability}, there exists an automorphism $\phi:L_p(L_q) \rightarrow L_p(L_q)$ such that $\phi \circ g_1 = g_2$.  Note then that $\| \phi(f_1(e_k)) - f_2(e_k) \| \leq \|\phi(f_1(e_k) - g_1(e'_k)) \| + \| f_2(e_k) - g_2(e'_k) \| < \varepsilon$.  
\end{proof}

In a manner similar to that of Theorem \ref{t:finite-generated-blplq-homogeneity}, we can also extend the AUH property to finitely generated sublattices of $L_p(L_q)$ as well: 

\begin{theorem}\label{t:finite-generated-homogeneity}
	For all $1\leq p, q<\infty$ with $p/q \notin \N$, The lattice $L_p(L_q) $ is AUH for the class $\overline{ \K{p}{q}}$ of its finitely generated lattices.
\end{theorem}

\begin{proof}
	Suppose $E \subseteq L_p(L_q)$ is finitely generated. Then since $E$ is order continuous and separable, it is the inductive limit of finite dimensional lattices as well, so pick a finite dimensional $E'$ with elements sufficiently approximating the generating elements of $E$, and proceed with the same proof as in Theorem \ref{t:finite-generated-blplq-homogeneity}.
\end{proof}

The argument used in Corollary \ref{c:homogeneity_disjointness} can also be used
to show:

\begin{corollary}\label{c:fg_homogeneous}
	
For $p/q \notin \N$,  $L_p(L_q)$ is disjointness preserving AUH over $\overline{\K{p}{q}}$.
\end{corollary}

\begin{remark}
	$L_p(L_q)$ for $p/q\notin \N$ is AUH over the entire class of its finitely generated sublattices, a property which is equivalent to such a class being a metric \textit{Fra\"iss\'e class} with $L_p(L_q)$ as its \textit{Fra\"iss\'e limit}. Recall that a class $\mathcal K$ of finitely generated lattices is \textit{Fra\"iss\'e} if it satisfies the following properties:
	\begin{enumerate}
		\item  \textit{Hereditary Property} (HP): $\mathcal K$ is closed under finitely generated sublattices.
		\item \textit{Joint Embedding Property} (JEP): any two lattices in $\mathcal K$ lattice embed into a third in $\mathcal K$. 
		
		\item  \textit{Continuity Property} (CP):  any lattice operation symbol are continuous with respect to the Fra\"iss\'e pseudo-metric $d^\mathcal K$ in \cite[Definition 2.11]{ben15}.
		
		\item  \textit{Near Amalgamation Property } (NAP): for any lattices $ E =< e_1,...e_n>_L$, $ F_1$ and  $ F_2$ in $\mathcal K$ with lattice embeddings $f_i:  E \rightarrow F_i$, and for all $\varepsilon > 0$, there exists a $ G \in \mathcal K$ and embeddings $g_i :  F_i \rightarrow  G$ such that $\|g_1\circ f_1( e_k) -  g_2 \circ f_2 ( e_k) \| < \varepsilon.$
		
		\item  \textit{Polish Property} (PP): The Fra\"iss\'e pseudo-metric $d^{\mathcal K}$ is separable and complete in $\mathcal K_n$ (the $\mathcal K$-structures generated by $n$ many elements). \\
	
	\end{enumerate}

	Now clearly the finitely generated sublattices of $L_p(L_q)$ fulfill the first two properties, and the third follows from the lattice and linear operations having moduli of continuity independent of lattice geometry.  In addition, if one can show that the class $\mathcal K$ has the $NAP$, has some separable $X$ which is universal for $\mathcal K$, and its NAP amalgamate lattices can be chosen so that they are closed under inductive limits, then one can prove that $\mathcal K$ also has the Polish Property (a technique demonstrated in \cite[Theorem 4.1]{tursi20} and more generally described in Section 2.5 of \cite{lup16}).  The main difficulty in proving that a class of lattices $\mathcal K$ is a Fra\"iss\'e class is in showing that it has the NAP.  However, thanks to Theorem \ref{t:finite-generated-homogeneity}, we have 
	
	\begin{corollary}
		$\overline{ \K{p}{q}}$ has the NAP.
	\end{corollary}
	
	Theorem \ref{t:finite-generated-homogeneity} implies an additional collection of AUH Banach lattices to the currently known AUH Banach lattices: namely $L_p$ for $1\leq p <\infty$,the Gurarij M-space $\mathcal M$  discovered in \cite{fer20}, and the Gurarij lattice discovered in \cite{tursi20}. \\
	
	However, if one considers classes of finite dimensional Banach spaces with Fra\"iss\'e limits using linear instead of lattice embeddings, the only known separable AUH Banach spaces are the Gurarij space and $L_p$ for $p \neq 4,6,8,...$, and it is currently unknown if there are other Banach spaces that are AUH over its finite dimensional subspaces with linear embeddings. Certain combinations of $p$ and $q$ are also ruled out for $L_p(L_q)$ as a potential AUH candidate as discussed in Problem 2.9 of \cite{fer20}: in particular, when $1\leq p,q < 2$, $L_p(L_q)$ cannot be linearly AUH.
\end{remark}

\section{Failure of homogeneity for $p/q\in \N$}\label{s:non_homgeneity}

Recall that when $E = <e_1,...,e_n> \in B\K{p}{q}$ is embedded into $L_p(L_q)$ through $f_1, f_2$, then we can achieve almost commutativity for any $p\neq q$.  However, the automorphism in Theorem \ref{t:blplq_finite_homogeneity} clearly preserves the equimeasurability of the generating basic atoms of $f_i(E)$ as it fixes $\mathbf 1$.\\

 In this section, we show that the results of Section \ref{s:equimeasurability} do not hold when $p/q \in \N$.  The first results in this section show that when some $e \in L_p(L_q)_+$ is sufficiently close to $\mathbf 1$, the automorphism originally used in the argument of Proposition \ref{p:homogeneity-one-element} sending $\mathbf 1$ to $e$ also perturbs selected functions piecewise continuous on their support in a controlled way.   Second, Theorem \ref{theorem:lindeExtension} does not hold, and thus we cannot infer equimeasurability for arbitrary finite dimensional sublattices of $L_p(L_q)$. Finally, we use these results to strengthen the homogeneity property for any $L_p(L_q)$ lattice assumed to be AUH, and then show that when $p/q \in \N$, $L_p(L_q)$ does not fulfill this stronger homogeneity property, and thus cannot be AUH.

\begin{lemma}\label{l:homogeneity_approximate_equimeasurability}  Let $ 1\leq p \neq q < \infty$, and let $<f_1,...,f_n> \subseteq L_p(L_q)$ be such that $\sum f_i = \mathbf 1$.  Suppose also that for a.e. $x$, $f_k(x,\cdot) = \mathbf 1_{[g_{k}(x),g_{k+1}(x)]}$ where each $g_k$ has finitely many discontinuities. Let $\varepsilon > 0$, and let $e \in S(L_p(L_q))_+$ fully support $L_p(L_q)$. Consider 
	\[ \phi(f)(x,y) = f \bigg(\widetilde {N[e]}(x)_p, \frac{\widetilde e_x(y)_q}{N^q[e](x)} \bigg) e(x,y) \]
	which is the lattice isometry defined in Proposition \ref{p:homogeneity-one-element} mapping $\mathbf 1$ to $e$.

	Then there exists $\delta$ such that if $\| \mathbf 1 - e\|<\delta$, then for each $k$, we have that $\|\phi(f_k) - f_k\| < \varepsilon$. 
	
\end{lemma}

\begin{proof}
	
	We can assume $\varepsilon < 1$. Let $K\subseteq [0,1]$ be a closed set such that for $1\leq k \leq n+1$, $g_k|_K$ is continuous and $\mu(K) > 1 - \varepsilon$. Pick $\delta'<\varepsilon$ such that for any $x,x' \in K$, if $|x-x'|< \delta' $, then $|g_k(x)-g_k(x')| < \varepsilon/4$. Now, let $\delta < {\delta'}^{2p}$ be such that $1-\frac{\delta'}{4} \leq (1-\delta)^p< (1+\delta)^p < 1+\frac{\delta'}{4}$, and suppose $\|\mathbf 1 - e\| < \delta$. Observe that for each $x$, we have $\widetilde{ N[\mathbf 1 - e]}(x)_p <\delta$. For each $1\leq k\leq n$, let 
	$$\widetilde f_n(x,y) = f\bigg(\widetilde{ N[e]}(x)_p, \frac{\widetilde e_x(y)_q}{N^q[e](x)} \bigg).$$
	Observe that $\|\widetilde f_k - \phi(f_k) \| < \delta < \varepsilon/4$, so it is enough to show that $\|\widetilde f_k - f_k\|$ is sufficiently small as well.  \\
	
	To this end, first note that since $f$ is being composed with increasing continuous functions in both arguments, each $\widetilde f_n(x, \cdot )$ is also the characteristic function of an interval: indeed, we have piecewise continuous $\widetilde g_1,..., \widetilde g_{n+1}$ with $\widetilde g_k(x) :=g(\widetilde{ N[e]}(x)_p)$ and $\widetilde g_{n+1}(x) = 1$ such that for each $k$, $\widetilde f_k(x,y) = \mathbf 1_{[\widetilde g_k(x), \widetilde g_{k+1}(x)]}(y)$. Also observe that for $M := \{x\in K: N[e-1](x) < \delta\}$, we have $\mu(M)) > 1- \delta' - \varepsilon$. In addition, as 
	
	$$ \|f_k -\widetilde f_k \|^p = \| N[f_k - \widetilde f_k]\|^p_p = \int \mu(D(x))^p \ dx,$$
	
	Where $D_k(x) = \{y: f_k(x,y) \neq \widetilde f_k(x,y) \}$.  The above set up, in combination with the triangle inequality properties of $N$, leads us to the following inequalities:
	
	\begin{itemize}
		\item For all $0\leq x\leq 1$, $|\widetilde {N[e]}(x)_p - x| < {\delta} $. 
		\item For all $x\in M$ , $|N[e](x) - 1| < \delta$.
		\item For all $x\in M$ and $0\leq y \leq 1$, $|\widetilde e_x(y)_q - y |< \frac{\delta'}2$.
		\item For all $x\in M$ and $0 \leq y \leq 1$, if $y':= \frac{\widetilde{e}_x(y)_q}{N^q[e](x)}$, then $| y' - e_x(y)_q| < \frac{\delta'}2$	(which implies with the above that $|y-y'| < \delta'$).
	\end{itemize}
	
	We now show that the above implies that $D_k(x) < 2\varepsilon$.  Observe first that for all $x\in M$, if $f_k(x,y) \neq \widetilde f_k(x,y)$ it must be because, but  $y' \notin [\widetilde g_k(x), \widetilde g_{k+1}(x)]$, or vice versa. In either case, it can be shown that either $|y-g_k(x)|< \delta +\frac{\varepsilon}{4}$ or $|y-g_{k+1}(x)|< \delta +\frac{\varepsilon}{4}$. Suppose $y\in[g_k(x), g_{k+1}(x)]$ and $y' < \widetilde g_k(x)$ (a similar proof will work in the case that $y' > \widetilde g_{k+1}(x)$. Then since $y >g_k(x)$,  $|y-y'|\leq \delta' $, and $|g_k(x) - \widetilde g_k(x)| < \frac{\varepsilon}{4}$, 
	 $$0 \leq y-g_k(x) = (y-y') +(y'-\widetilde g_k(x)) + (\widetilde g_k(x) - g_k(x)) < \delta +\frac{\varepsilon}{4}. $$
	 
	 It follows then that accounting for both ends of the interval $[g_k(x), g_{k+1}(x)]$ and for $x\in M$, we have $D_k(x) < 2\varepsilon$.  Resultantly,
	\begin{align*}
		\| f_k - \widetilde f_k \|^p = \int_M \mu(D(x))^p \ dx + \int_{ M^c}\mu(D(x))^p \ dx < (2\varepsilon)^p + \delta^p < 3\varepsilon^p,
	\end{align*}
	
	which can be made arbitrarily small.
	
\end{proof}

\begin{theorem}\label{t:homogeneity_approximate_equimeasurability}
	Let $1\leq p \neq q <\infty$ and suppose $L_p(L_q)$ is AUH over its finite dimensional sublattices.  Let $f_i:E\rightarrow L_p(L_q)$ be lattice embeddings with $E = <e_1,...,e_n>$ such that $f_i(x) = \mathbf 1$ for some $x\in E$.  Then for all $\varepsilon > 0$, there exists an automorphism $\phi$ fixing $\mathbf 1$ such that $\| \phi f_1 - f_2 \| < \varepsilon$.
\end{theorem}

\begin{proof}
 Assume the above, and pick $E' = <e'_1,...,e'_m> \subseteq L_p(L_q)$, where $e'_k = a_k \cdot \mathbf 1_{A_k\times B_k}$ with $A_k$ and $B_k$ intervals such that $\sum_k \mathbf 1_{A_k\times B_k} = \mathbf 1$ and for each $e_k$ there is $x_k \in S(E')_+$ such that $\|x_k - f_2(e_k)\| < \frac{\varepsilon}{4n}$. \\

Since $L_p(L_q)$ is AUH, there exists an automorphism $\psi$ such that $\| \psi f_1 - f_2 \| < \delta$, where $\delta$ satisfies the conditions for $\frac{\varepsilon}{4mn}$ and each of the $e'_k$'s in $E'$ in Lemma \ref{l:homogeneity_approximate_equimeasurability}.  Now pick the automorphism $\phi'$ over $L_p(L_q)$ mapping $\mathbf 1$ to $\psi f_1(x)$ as defined in Lemma \ref{l:homogeneity_approximate_equimeasurability}.  It follows that for each $e'_k$, $\| \phi'(e'_k) - e'_k \| < \frac{\varepsilon}{4mn}$, so $\|\phi'(x_k) - x_k \| < \frac{\varepsilon}{4n}$.  Thus for each $e_k \in E$,

\begin{align*} \| \phi' f_2(e_k) - \psi f_1(e_k) \|  \leq & \| \phi' (f_2(e_k) -x_k) \| +\| \phi'(x_k) - x_k \| \\ 
	 + & \| x_k - f_2(e_k) \| + \| f_2(e_k)  - \psi f_1(e_k) \| < \frac{\varepsilon}{n}, 
\end{align*}

Now let $\phi = {\phi'}^{-1}\circ \psi$ to obtain the desired automorphism; then $\|\phi f_1 - f_2\| < \varepsilon.$

\end{proof}

	The above can be used to show that if $L_p(L_q)$ is AUH and $f_i(E)$ contains $\mathbf 1$ for $i=1,2$, then we can induce almost commutativity with automorphisms fixing $\mathbf 1$ as well.  This will allow us to reduce possible automorphisms over $L_p(L_q)$ to those that in particular fix $\mathbf 1$. The importance of this result is that these particular homomorphisms fixing $\mathbf 1$ must always preserve base-equimeasurability for characteristic functions, as shown in Proposition \ref{p:1-1-isometry_representation}.  Thus a natural approach in disproving that $L_p(L_q)$ is AUH would involve finding sublattices containing $\mathbf 1$ which are lattice isometric but whose generating elements are not base-equimeasurable.  The following results do exactly that:

\begin{lemma}\label{l:p/q_exception}
	Lemma \ref{lemma:LindeExtention} fails when $r:=p/q\in \N$. In particular, there exists a non-zero measure $\nu:=\alpha - \beta$, with $\alpha$ and $\beta$ positive measures such that for all polynomials $P$ of degree $j \leq r$, 
	
	$$ \int_0^1 P(x) \ d\nu(x) = 0.$$ 
\end{lemma}

\begin{remark}
It is already known that a counter-example exists for $L_r(0,\infty)$ for all $r\in \N$, with 
$$d\nu(x) = e^{-u^{\frac14}} \sin (u^{\frac14}) \ du$$
(see \cite{rudin1976p} and \cite{Linde1986UniquenessTF} for more details).
\end{remark}.  

Here we provide another example over the unit interval:

\begin{proof}
	 Fix such an $r$, and define a polynomial $g(x)$ of degree $r+1$ with $g(x) = \sum_0^{r+1} a_i x^i$ such that for all $0 \leq j \leq r$, $\int_0^1 x^j g(x) \  dx = 0$.  This can be done by finding a non-trivial $a_0,..., a_{r+1}$ in the null set of the $(n+1) \times (n+2)$ size matrix $A$ with $A(i,j) = \frac{1}{i+j+1}$.  Then let $d\nu(x) = g(x) \ dx$. Let $\alpha = \nu_+$ and $\beta = \nu_-$.  Clearly $\alpha$ and $\beta$ are finite positive Borel measures, but since $g \neq 0$, $\alpha \neq \beta$.
\end{proof}

\begin{lemma}\label{l:p/q_nonequimeasurable}
	Let $p/q \in \N$.  Then there exists a two dimensional lattice $E = <e_1, e_2>$ and lattice embeddings $f_i:E\rightarrow L_p(L_q)$ with $ \mathbf 1 \in E$ such that $g_1(\mathbf {e})$ and $g_2(\mathbf{e})$ are not base-equimeasurable.
\end{lemma}

\begin{proof}
	 Let $f(x)$ be a polynomial of degree at least $r+1$ as defined in Lemma \ref{l:p/q_exception} such that for all $0 \leq k \leq r$, $\int_0^1 t^k f(t) \ dt = 0$, and $\int_0^1 |f(x)| \ dx = 1$.   Let $h_1(x) = \frac12 + f(x)_+$, and let $h_2(x) = \frac12 +f(x)_-$.  Note that each $h_i(x) > 0$, and furthermore that $\int_0^1 h_i(t) \ dt = 1$.  Additionally, each map $H_i(x) = \int_0^x h_i(t) \ dt$ is strictly increasing with $H_i(0) = 0$ and $H_i(1) = 1$.  Now we will construct characteristic functions $f^i_j \in L_p(L_q)$ such that the linear map $f_j^1 \mapsto f_j^2$ induces an isometry, but $N\mathbf f^1$ and $\mathbf f^2$ are not base-equimeasurable.  From there, we let $e_j = \frac{f^1_j}{\|f^i_j\|}$, and let $g_i$ be the lattice isometry induced by $g_i(e_j) = \frac{f^i_j}{\|f^i_j\|}$, \\
	 
	 To this end, let
	\[F^i_1(x) := H_i^{-1}(x), \text{ and } F^i_2(x) := 1 - F^i_1(x).  \]
	
	  Observe that $F^1_1(x) \neq F^2_1(x)$. Indeed, one can show that the associated push forwards $d {F^i_1} _\# \mu$ for each $F^i_1$ have the corresponding equivalence:
	\[ d {F^i_1} _\# \mu (x)  = h_i(x) \ dx   \]
	So $(F_1^1, F_2^1)$ and $(F_1^2, F_2^2)$ are not equimeasurable.  However, For $0 \leq j \leq r$, $u^j h_i(u) \ du = u^j \ d{F_1^i}_\#(u)  = F_1^i(x) ^j \ dx$, so it follows from the construction of the $h_i$'s that 
	$$ \int_0^1 F_1^1(x) ^j \ dx = \int_0^1 F_1^2(x)^j \ dx.$$
	Thus for any $v_1, v_2 > 0$, since $F^i_1$ and $F^i_2$ are both positive, we have 
	
	\begin{align*}
		& \int_0^1 |v_1 F^1_1(x) + v_2 F^1_2(x) |^r \ dx =  \int_0^1 ((v_1 - v_2) F^1_1(x) + v_2)^r \ dx \\ 
		= & \sum_0^r \binom{r }{j} (v_1 - v_2)^jv_2^{r-j} \int_0^1 F^1_1(x)^j \ dx 	=  \int_0^1 |v_1 F^2_1(x) + v_2 F^2_2(x) |^r \ dx
	\end{align*}
	
	To conclude the proof, let $f^i_1(x,y) = \mathbf 1_{[0, F^i_1(x)]}(y)$, and let $f^i_2 = \mathbf 1 - f^i_1$.  Clearly $N[f^i_j] = F^i_j$.
\end{proof}

\begin{theorem}\label{t:not-AUH}
	If $p/q \in \N$ and $p\neq q$, then $L_p(L_q)$ is not AUH for the class of its finite dimensional sublattices.
	\end{theorem}

\begin{proof}
	Fix $p/q\in \N$, and let $E$ be the 2-dimensional lattice generated in Lemma \ref{l:p/q_nonequimeasurable}, with $f_i:E\rightarrow L_p(L_q)$ embeddings mapping to copies of $E = <e_1,e_2>$ such that $f_1(\mathbf e)$  and $f_2(\mathbf e)$ are not base-equimeasurable. In addition, by assumption $\mathbf 1 \in E$.  For notational ease, let $F^i_j = N[f_i(e_j)]$.  \\
	
	Suppose for the sake of contradiction that $L_p(L_q)$ is AUH. Pick some measurable $C \subseteq [0,1]^2$ and $\varepsilon > 0$ such that $$*  \quad \mathbf F^2_\# \mu (C)   > \mathbf F^1_\# \mu(C+\varepsilon)+ \varepsilon, $$
	
	where $$C+\varepsilon = \{ \mathbf t \in [0,1]^2 : \|\mathbf {t-s}\|_\infty < \varepsilon \text{ for some }\mathbf s \in C \}.$$
	 By Theorem \ref{t:homogeneity_approximate_equimeasurability}, there is some lattice automorphism $\phi:L_p(L_q)\rightarrow L_p(L_q)$ fixing $\mathbf 1$ such that $\| \phi\circ f_1 - f_2\| < \varepsilon^2$. Let $ \phi F^i_j = N[\phi f_i(e_j)] $.  By Proposition \ref{p:1-1-isometry_representation}, $\phi$ preserves base-equimeasurability, so for any measurable $B$, 
$$\quad \phi \mathbf F^1_\# \mu (B)  = \mathbf F^1_\# \mu(B). $$
By the properties of $N$, we also have $\|\phi F^1_j  - F^2_j\|_p \leq \|\phi f_1(e_j) - f_2(e_j) \|$. It also follows that 
$$\mu (t: \| \phi \mathbf F^1(t) - \mathbf F^2(t) \|_\infty > \varepsilon ) < \varepsilon, $$
 so $\phi \mathbf F^1_\# \mu (C+\varepsilon) +\varepsilon > \mathbf F^2_\# \mu(C)$, but this contradicts the assumption (*).  So Theorem \ref{t:homogeneity_approximate_equimeasurability} cannot apply, implying that $L_p(L_q)$ is not AUH as desired.
\end{proof}

\begin{remark}
	For $p/q\in \N$, $L_p(L_q)$ is the unique lattice that is separably AUH over finitely generated $BL_pL_q$ spaces, since up to isometry it is the unique doubly atomless $BL_pL_q$ space.  In light of Theorem \ref{t:not-AUH}, this implies that the class of finitely generated sublattices of $L_p(L_q)$ is not a Fra\"iss\'e class as defined in \cite{ben15}, as $L_p(L_q)$ is the only possible candidate as a Fra\"iss\'e limit. \\
	
	In particular, $L_pL_q$ lacks the NAP.  Indeed, otherwise, one can use that NAP with $BL_pL_q$ amalgamate lattices and \cite[Proposition 2.8]{henson11} to situate a $d^\mathcal K$-Cauchy sequence into a Cauchy-sequence of generating elements in an ambient separable $BL_pL_q$ lattice. Thus $\overline{ \K{p}{q}}$ would also have the Polish 
	Property, implying that $\overline{\K{p}{q}}$ is a Fra\"iss\'e class. Since the only possible candidate Fra\"iss\'e limit space is $L_p(L_q)$ itself, this would contradict Theorem \ref{t:not-AUH}.
\end{remark}

\bibliographystyle{abbrv}
\bibliography{bibliography}
\end{document}